\documentclass{amsart}
\usepackage{epsfig}

\newtheorem{theorem}{Theorem}[section]
\newtheorem{lemma}[theorem]{Lemma}
\newtheorem{proposition}[theorem]{Proposition}
\theoremstyle{definition}

\theoremstyle{remark}

\numberwithin{equation}{section}

\def\sv{\mathcal V}
\def\barv{\bar{v}}
\def\bars{\bar{s}}

\begin{document}

\title{Rigidity of two-dimensional Coxeter groups}

\author[Patrick Bahls]{Patrick Bahls}

\address{Department of Mathematics \\ University of Illinois at 
Urbana-Champaign \\ Urbana, IL 61801}

\email{pbahls@math.uiuc.edu}



\begin{abstract}
A Coxeter system $(W,S)$ is called {\it two-dimensional} if the Davis 
complex associated to $(W,S)$ is two-dimensional (equivalently, every 
spherical subgroup has rank less than or equal to 2).  We prove that 
given a two-dimensional system $(W,S)$ and any other system $(W,S')$ 
which yields the same reflections, the diagrams corresponding to these 
systems are isomorphic, up to the operation of diagram twisting defined by 
Brady, McCammond, M\"uhlherr, and Neumann in \cite{BMMN}.  As a step in 
the proof of this result, certain two-dimensional groups are shown to be 
reflection rigid in the sense of \cite{BMMN}, and a result concerning the 
strong rigidity of two dimensional systems is given in the final section.
\end{abstract}

\thanks{The author was supported by an NSF VIGRE postdoctoral grant.}

\maketitle

\section{Introduction}

A {\it Coxeter system} is a pair $(W,S)$ where $W$ is a group with a 
presentation of the form $\langle S \ | \ R \rangle$, $S=\{s_i\}_{i \in 
I}$, and

$$R = \{ (s_is_j)^{m_{ij}} | m_{ij} \in \{1,2,...,\infty\}, m_{ij}=m_{ji}, 
\ {\rm and} \ m_{ij}=1 \Leftrightarrow i=j \}.$$

\noindent When $m_{ij} = \infty$, the element $s_is_j$ has infinite 
order.  A group $W$ with such a presentation is called a {\it Coxeter 
group}, and $S$ is called a {\it fundamental generating set}.

Let $T \subseteq S$.  Denote by $W_T$ the subgroup of $W$ generated by the 
elements in $T$.  Such a subgroup is called a {\it standard parabolic 
subgroup} of $W$, and any conjugate of such a group is called a {\it 
parabolic subgroup}.  If $W_T$ is finite, $W_T$ is called a {\it spherical 
subgroup}.  It is well-known (see \cite{Bo}, for instance) that $(W_T,T)$ 
is a Coxeter system for any subset $T \subseteq S$, and therefore $W_T$ is 
a Coxeter group in its own right, with the obvious presentation.  It is 
also known that any spherical subgroup $W_T$ contains a unique longest 
element with respect to the set $S$ (see \cite{Bo}), which we denote by 
$\Delta_T$.  This element has the property that $\Delta_T$ conjugates any 
element $t \in T$ to some $t' \in T$.

The information contained in the presentation $\langle S \ | \ R \rangle$ 
above can be displayed nicely by means of a {\it Coxeter diagram}.  The 
Coxeter diagram $\sv$ associated to the Coxeter system $(W,S)$ is an 
edge-labeled graph whose vertices are in one-to-one correspondence with 
the generating set $S$ and for which there is an edge $[s_is_j]$ labeled 
$m_{ij}$ between two vertices $s_i$ and $s_j$ if and only if $i \neq j$ 
and $m_{ij} < \infty$.

Given a spherical subgroup $W_T$ of $S$, it is clear that the subgraph of 
$\sv$ induced by the generators in $T$ is a simplex in the combinatorial 
sense.  We call such a simplex a {\it spherical simplex}, and say that it 
is {\it maximal} if it is not properly contained in another spherical 
simplex.

In the sequel, we frequently omit the word ``Coxeter'' when discussing 
groups, systems, and diagrams, as these words will be used in no other 
context.

It is easy to see that the diagram fully and faithfully records all of the 
information in the presentation $\langle S \ | \ R \rangle$.  It is also 
easy to see that to a given group $W$ there may correspond more than one 
system (and therefore diagram).  For example, the dihedral group $D_{2k}$ 
of order $4k$ has the presentations

$$\langle a,b \ | \ a^2, b^2, (ab)^{2k} \ \rangle$$

\noindent and

$$\langle c,d,g \ | c^2, d^2, g^2, (cd)^2, (cg)^2, (dg)^k \rangle$$

\noindent when $k$ is odd.  These correspond to diagrams consisting of a 
single edge labeled $2k$, and a triangle with edge labels $\{2,2,k\}$, 
respectively.

Therefore one may consider the question: to what extent is a given Coxeter 
system unique?  As a first step toward answering this question, we must 
decide what is meant by ``unique''.

We say that the group $W$ is {\it rigid} if given any two systems $(W,S)$ 
and $(W,S')$, there is an automorphism $\alpha \in {\rm Aut}(W)$ 
satisfying $\alpha(S)=S'$.  Equivalently, the diagrams corresponding to 
these two systems are isomorphic as edge-labeled graphs.  We say that $W$ 
is {\it strongly rigid} if such an automorphism $\alpha$ can always be 
chosen to lie in ${\rm Inn}(W)$; i.e., any two fundamental generating sets 
are conjugate to one another.

We can relax these conditions slightly.  We require the notion of a {\it 
reflection}.  A reflection in the system $(W,S)$ is any conjugate 
$wsw^{-1}$ of a generator $s \in S$.  We say that a Coxeter system $(W,S)$ 
is {\it reflection rigid} if given any other system $(W,S')$ which yields 
the same reflections, there is an automorphism $\alpha$ of $W$ satisfying 
$\alpha(S)=S'$.  Finally, $(W,S)$ is said to be {\it strongly reflection 
rigid} if given any other system $(W,S')$ yielding the same reflections, 
such an automorphism $\alpha$ can be found in ${\rm Inn}(W)$.  We call $W$ 
{\it reflection independent} if every two systems for $W$ yield the same 
reflections.  Clearly if $W$ is reflection independent, then (strong) 
rigidity and (strong) reflection rigidity are equivalent.

A number of results have been proven that characterize the groups 
that satisfy these rigidity conditions.  Furthermore, there are other 
characterizations of uniqueness with which we will not concern ourselves 
in this paper.  (See \cite{Ba1}, \cite{Ba2}, \cite{Ba3}, \cite{BaMi1}, 
\cite{BMMN}, \cite{ChDa}, \cite{Ka}, \cite{Mi}, \cite{MuWe}, \cite{Ra}.)

In this paper we will generalize the method used in \cite{Ba3} in order to 
describe the extent to which {\it two-dimensional} Coxeter groups are 
rigid.  A system $(W,S)$ is called two-dimensional (or 2-d) if no three 
distinct generators from $S$ generate a finite subgroup of $W$.  (The term 
``two-dimensional'' refers to the dimension of the Davis complex, a 
simplicial complex associated to the system $(W,S)$.  See \cite{ChDa}, 
\cite{Da} for more details regarding this complex and its usefulness.)  
The group $W$ is called two-dimensional if there exists a two-dimensional 
system $(W,S)$.  (As a consequence of the main theorem below, we will see 
that this distinction is unnecessary in the presence of reflection 
independence.)  In order to describe the results we obtain, we must 
introduce the important notion of {\it diagram twisting}, due to Brady, 
McCammond, M\"uhlherr, and Neumann (in \cite{BMMN}).

Given a Coxeter system $(W,S)$, suppose that $T$ and $U$ are disjoint 
subsets of $S$ satisfying

\vskip 2mm

\noindent 1. \ $W_U$ is spherical, and

\vskip 2mm

\noindent 2. \ every vertex in $S \setminus (T \cup U)$ which is connected 
to a vertex of $T$ by an edge is also connected to every vertex in $U$, by 
an edge labeled 2.

\vskip 2mm

Under these conditions, we may define a new diagram (and therefore new 
system) $\sv'$ for $W$ by changing every edge $[tu]$ ($t \in T$, $u \in 
U$) to an edge $[tu']$, where $u' = \Delta_U^{-1} u \Delta_U$, leaving 
every other edge unchanged.  This modification results in a generating set 
$S'$ obtained from $S$ by replacing $t \in T$ with $\Delta_U^{-1} t 
\Delta_U$.

This operation is called a {\it diagram twist}, because of the way that 
we ``twist'' around the subdiagram representing the group $W_U$.

We require a few new terms in order to state this paper's main results, 
remaining consistent with the terminology of \cite{MuWe}.  If $\sv$ is 
connected and $s$ is a vertex in $\sv$ such that $\sv \setminus \{s\}$ is 
disconnected, $s$ is called a {\it cut vertex} of $\sv$.  If $\sv$ has no 
such vertices, we say that $\sv$ is {\it one-connected}.  If $\sv$ is 
one-connected and there exists no edge $[st]$ such that $\sv \setminus 
[st]$ is disconnected, then $\sv$ is called {\it edge-connected}.  If 
$\sv$ is one-connected and there exists no edge $[st]$ with odd label such 
that $\sv \setminus [st]$ is disconnected, we call $\sv$ {\it 
odd-edge-connected}.  (Thus $\sv$ is odd-edge-connected if it is 
edge-connected.)

\begin{theorem} \label{main}
Let $(W,S)$ be a two-dimensional Coxeter system with diagram $\sv$.  Then 
$(W,S)$ is reflection rigid, up to diagram twisting.  (That is, given a 
system $(W,S')$ which yields the same reflections as $(W,S)$, there is a 
sequence of diagram twists which transforms the first system into the 
second.)
\end{theorem}

As a step in the proof of the main theorem, we will prove

\begin{theorem} \label{edgeconnected}
Let $(W,S)$ be a two-dimensional Coxeter system with odd-edge-connected 
diagram $\sv$.  Then $(W,S)$ is reflection rigid.
\end{theorem}

Furthermore, we will prove a theorem (Theorem~\ref{sr}) concerning the 
strong rigidity of 2-d Coxeter groups.  Its statement will be deferred 
until the final section of this paper.

The above results partially generalize the similar results obtained by 
M\"uhlherr and Weidmann in \cite{MuWe}.  Indeed, in this paper we will 
make similar use of the results of \cite{MiTs} in order to complete the 
proof of Theorem~\ref{main} (see \cite{MuWe}, Section 8).  However, the 
preliminary arguments are very different, and will be introduced in the 
following section.  The author has also recently learned that T. Hosaka 
has proven independently a slightly weaker result concerning rigidity of 
two-dimensional Coxeter groups.

This is the third paper in a series (see \cite{Ba3}, \cite{Ba4}) which 
makes use of similar techniques in order to establish structural 
properties of Coxeter groups.  It is clear that these techniques can be 
pushed even further to prove results about yet more general Coxeter 
groups.  This will be done in subsequent papers.

The author gratefully acknowledges helpful discussions with Ruth Charney, 
Ilya Kapovich, and Richard Weidmann during the writing of this paper.

\section{Circuits and centralizers}

We begin by sketching the argument that we will use to prove 
Theorem~\ref{main}.  Let $(W,S)$ be a 2-d system, and let $(W,S')$ be 
another system for $W$.  Denote the corresponding diagrams by $\sv$ and 
$\sv'$.  Let us assume until further mention that $W$ is reflection 
independent.  Our goal is to show that, up to twisting, $\sv$ and $\sv'$ 
are identical.

The two-dimensionality of $(W,S)$ allows us to establish a matching 
between the edges of $\sv$ and the edges of $\sv'$, using the following 
result from \cite{Ka}.

\begin{theorem} \label{ematch}
Let $W$ be a Coxeter group with diagram $\sv$, all of whose maximal 
spherical simplices are of the same dimension.  Then, given any other 
Coxeter system $(W,S')$ with diagram $\sv'$, there is a one-to-one 
correspondence $\phi$ between the maximal spherical simplices of $\sv$ and 
those of $\sv'$.  Moreover, for any maximal spherical simplex $\sigma$ in 
$\sv$, there is an element $w \in W$ such that $wW_{\sigma}w^{-1} = 
W_{\phi(\sigma)}$.
\end{theorem}

In our case, every maximal spherical simplex is an edge, and therefore has 
dimension 1.  We apply Theorem~\ref{ematch} to obtain a matching between 
the edges of $\sv$ and the edges of $\sv'$ which respects conjugacy as 
indicated in the theorem.

Why must each edge of $\sv$ be matched with an {\it edge} of $\sv'$?  If 
there were an edge $[st]$ in $\sv$ such that $wW_{[st]}w^{-1}=W_{\sigma}$ 
for some $\sigma$ of dimension $>1$, then $D_n \cong W_{\sigma}$, where 
$n$ is the order of $st$.  However, it is an easy matter (see \cite{Ba1}) 
to show that this can only happen if $n=2k$, $k$ odd, and $\sigma$ is a 
triangle with edge labels $\{2,2,k\}$.  In this case, the central element 
of $W_{[st]}$ (which is of even length with respect to $(W,S)$, and is 
therefore not a reflection) is a reflection in $(W,S')$, contradicting the 
assumption that $W$ is reflection independent.  An immediate corollary is 
that every system corresponding to $W$ is 2-d, so it matters not whether 
we refer to the group or to the system as 2-d, provided $W$ is reflection 
independent.  (Alternately, one may apply Lemma~1.5 of \cite{ChDa}.)

Let $\phi$ be the matching whose existence is guaranteed by 
Theorem~\ref{ematch}.

We now consider circuits in the diagram $\sv$.  A {\it simple circuit} of 
length $k$ in $\sv$ is a collection $C$ of $k$ distinct edges 
$\{[s_1s_2],...,[s_ks_1]\}$ for which $s_i \neq s_j$ when $i \neq j$.  
Define $d(i,j)={\rm min}\{|i-j|,k-|i-j|\}$ for $1 \leq i,j \leq k$.  We 
call a simple circuit $C$ {\it achordal} if for any two vertices $s_i 
\neq s_j$ in $C$ such that $d(i,j)>1$, $[s_is_j]$ is not an edge in $\sv$.

We shall prove the following theorem.

\begin{theorem} \label{cmatch}
Let $(W,S)$, $(W,S')$, $\sv$, $\sv'$, and $\phi$ be as above.  Let $C$ be 
an achordal circuit of length $k$ in $\sv$, as above.  Then there is an 
achordal circuit $C' = \{[\hat{s}_1\hat{s}_2],...,[\hat{s}_k\hat{s}_1]\}$ 
in $\sv'$ such that $\{\hat{s}_i,\hat{s}_{i+1}\} = \phi(\{s_i,s_{i+1}\})$ 
for $1 \leq i \leq k$.  Moreover, for each edge $[s_is_{i+1}]$ there is an 
element $w_{i+1} \in W$ such that $w_{i+1}s_iw_{i+1}^{-1}=\hat{s}_i$ and 
$w_{i+1}s_{i+1}w_{i+1}^{-1}=\hat{s}_{i+1}$ both hold.
\end{theorem}

Therefore not only do the edges match up nicely, but the achordal circuits 
do as well.  In fact, we can do better:

\begin{theorem} \label{cmatch2}
Let $(W,S)$, $(W,S')$, $\sv$, $\sv'$, $\phi$, $C$, and $C'$ be as in 
Theorem~\ref{cmatch}.  Let $s_i$ and $s_j$ be distinct vertices on $C$, 
with $\{e_i=[s_is_{i+1}],[s_{i+1}s_{i+2}],...,e_j=[s_{j-1}s_j]\}$ a 
subpath of $C$ between them.  Let $w_i$ and $w_j$ be the group 
elements which conjugate the edges $e_i$ and $e_j$, respectively, to their 
corresponding edges in $C'$.  Then $w_iw_j^{-1}$ can be written 
$\alpha_1\alpha_2 \cdots \alpha_r$, where for every $l$ ($1 \leq l \leq 
r$) $\alpha_l \in S'$ and one of the following holds.

\vskip 2mm

\noindent 1. \ For every $l'$ ($1 \leq l' \leq k$), $\alpha_l \neq 
\hat{s}_{l'}$, and $\alpha_l$ commutes with at least 2 vertices which lie 
on $C'$.  Moreover, we can find two such elements, $\hat{s}_{l_1}$ and 
$\hat{s}_{l_2}$, such that the path $\{ [\hat{s}_{l_1}\alpha_l], 
[\alpha_l\hat{s}_{l_2}]\}$ separates $C'$ into two circuits, one 
containing $\hat{s}_i$ and the other containing $\hat{s}_j$.

\vskip 2mm

\noindent 2. \ $\alpha_l=\hat{s}_i$, in which case both 
$[\hat{s}_{i-1}\hat{s}_i]$ and $[\hat{s}_i\hat{s}_{i-1}]$ are labeled 2, 
or $\alpha_l=\hat{s}_j$, in which case both $[\hat{s}_{j-1}\hat{s}_j]$ and 
$[\hat{s}_j\hat{s}_{j+1}]$ are labeled 2.
\end{theorem}

Although Theorem~\ref{cmatch2} appears very technical, it addresses 
precisely the issues that must be faced when dealing with strong rigidity 
in the presence of edges labeled 2.  (Compare the arguments of Section 6 
in \cite{Ba3}; in particular, those used in Cases 1 and 2.)  We note that 
in case no edges in $\sv$ are labeled 2, $w_i=w_j$ must hold for all edges 
$w_i$ and $w_j$; thus the circuit $C$ is in this case ``strongly rigid''.

M\"uhlherr and Weidmann also consider achordal circuits in \cite{MuWe}, 
but their approach to these circuits is very different from that adopted 
here, where we draw upon the techniques developed in \cite{Ba3} and 
\cite{Ba4}.

Once Theorem~\ref{cmatch} and Theorem~\ref{cmatch2} have been established, 
it will be a relatively straightforward matter to reconstruct the unique 
(up to twisting) diagram $\sv$ which is built up from the achordal 
circuits.

As will become clear, our analysis of the achordal circuits in $\sv$ 
will depend upon an understanding of the centralizer $C(s)$ of an 
arbitrary generator $s \in S$.  To that end, we recall in the next 
theorem the structure of $C(s)$ (first given in \cite{Br}).  We 
also introduce notation which will remain fixed throughout the remainder 
of the paper.

Let $(W,S)$ be an {\it arbitrary} Coxeter system and suppose $s,t \in S$ 
are elements of the fundamental generating set $S$.  If $m_{st} = 2k$ is 
even, denote by $u_{st}$ the element $(st)^{k-1}s$.  We note that $u_{st}$ 
commutes with $t$ (in fact, $u_{st}=s$ if $st=ts$).  If $m_{st} = 2k+1$ is 
odd, denote by $v_{st}$ the element $(st)^k$.  Note that $v_{st} s 
v_{st}^{-1} = t$.  More generally, there is an path in the diagram $\sv$ 
between two vertices $s$ and $t$ which consists entirely of odd edges 
if and only if $s$ and $t$ are conjugate to one another.  In fact, if $\{ 
[ss_1],[s_1s_2],...,[s_kt]\}$ is such a path, then $v_{s_kt} 
v_{s_{k-1}s_k} \cdots v_{ss_1}$ conjugates $s$ to $t$.

Let $\bar{\sv}$ be the graph resulting from a diagram $\sv$ by removing 
all edges with even labels.  As in \cite{Br}, we can identify elements of 
the fundamental group of $\bar{\sv}$ with paths in $\bar{\sv}$ which start 
and end at a fixed vertex $s \in S$ and which never backtrack.   For the 
fixed vertex $s \in S$, let ${\mathcal B}(s)$ be a collection of simple 
circuits in $\bar{\sv}$ containing $s$ such that ${\mathcal B}(s)$ 
generates the fundmental group of $\bar{\sv}$.  

The following was first proven by Brink in \cite{Br}.  (The generators 
given here can be computed by arguments similar to those in \cite{BaMi2}.)

\begin{theorem} \label{brink}
Let $(W,S)$ be an arbitrary Coxeter system with diagram $\sv$, and let $s 
\in S$.  Then $C(s)$ is the subgroup of $W$ generated by

$$\{s\} \cup A \cup B$$

\noindent where

$$A = \{ vu_{ts_k}v^{-1} \ | \ v = v_{s_1s} v_{s_2s_1} \cdots 
v_{s_ks_{k-1}}, t,s_i \in S, m_{ts_k} \ {\rm even}; m_{s_1s}, 
m_{s_is_{i-1}} \ {\rm odd} \}$$

\noindent and

$$B = \{ v_{s_1s} v_{s_2s_1} \cdots v_{ss_k} \ | \ 
\{[ss_1],...,[s_ks]\} \in {\mathcal B}(s) \}.$$
\end{theorem}

We will use this description of the centralizer $C(s)$ in the sequel.  

\vskip 1mm

\noindent {\bf Remarks.} \ When $(W,S)$ is 2-d, it can be shown that 
distinct choices of $s_1,s_2,...,s_k$ and $t$ in $A$ and $B$ above yield 
distinct generators.  This may not be the case if $(W,S)$ is not 2-d.  

Moreover, it is not difficult to compute a geodesic form for an element 
$w$ in $C(s)$.  To do this, first express $w$ as a product in the given 
generators.  Factor all occurrences of $s$ as a single letter to the end 
of of the word and cancel, yielding either $s$ or $1$.  Next, perform all 
``obvious'' cancellation; that is, given $\alpha=v_{s_1s} \cdots v_{ss_k}$ 
and $\beta=v_{s'_1s} \cdots v_{ss'_{k'}}$ in $B$, for some $i$, $0 \leq i 
\leq {\rm min}\{k,k'\}$ we have 
$s_k=s'_1,s_{k-1}=s'_2,...,s_{k-i+1}=s'_i$, so that $\alpha \cdot 
\beta = v_{s_1s} \cdots v_{s_{k-i+1}s_{k-1}} v_{s'_{i+1}s'_i} \cdots 
v_{ss'_{k'}}$.  Similar cancellation occurs in a product of two generators 
from $A$, and in a product of a generator from $A$ and a generator from 
$B$.

We claim that the word that results after such cancellation is geodesic.  
In order to prove this, we appeal to a result of Tits.  From Section 2 of 
\cite{Ti} we conclude that if the word resulting from the previous 
paragraph were not geodesic, we would be able to shorten the word by 
successively replacing subwords $(st)^n$ with $(ts)^n$ when $st$ has order 
$2n$ and subwords $(st)^ns$ with $(ts)^nt$ when $st$ has order $2n+1$, and 
then canceling any adjacent occurrences of the same letter which might 
arise in the course of these replacements.  However, thanks to 
two-dimensionality, no such shortening replacements can be performed, 
(perhaps) aside from commuting the single occurrence of $s$ that may occur 
at the end.

The above argument (replacing the one half of a relator with the other 
half) will be used again in the following sections.  We refer to the 
process of shortening a word $w$ in the manner described above as the {\it 
Tits process} (TP).

\section{Matching edges in a given circuit}

In this section we retrace the arguments from \cite{Ba3}, adapting them as 
necessary to the case of 2-d systems.  In fact, many of the arguments 
throughout the remainder of the paper will parallel arguments from 
\cite{Ba3} (such analogous arguments will be indicated).

Let $(W,S)$ be a 2-d system with diagram $\sv$, and let $(W,S')$ be 
another system, with diagram $\sv'$, yielding the same reflections as 
$(W,S)$.  We fix all of this notation for the remainder of the paper.

Let $C = \{[s_1s_2],...,[s_ks_1]\}$ be an achordal circuit in $\sv$.  By
Theorem~\ref{ematch}, for every $i=1,...,k$ there exists an edge 
$[s''_{i-i}s'_i]$ in $\sv'$ and an element $w_i \in W$ such that $w_i 
W_{[s_{i-1},s_i]} w_i^{-1} = W_{[s''_{i-1},s'_i]}$.  By considering the 
possible generators for the dihedral group $W_{[s''_{i-1},s'_i]}$, we can 
assume that

$$w_i s_{i-1} w_i^{-1} = s''_{i-1} \ {\rm and} \ w_i s_i w_i^{-1} = 
z_{i-1,i} s'_i z_{i-1,i}^{-1},$$

\noindent for some word $z_{i-1,i} \in W_{[s''_{i-1},s'_i]}$.  Let 
$m=m_{i-1,i}$ be the order of $s_{i-1}s_i$.  One may prove by direct 
computation that after suitably modifying $w_i$ we can assume
$z_{i-1,i} = (s'_is''_{i-1})^j$, where $0 \leq j \leq [\frac{m-1}{4}]$ if 
$m$ is even and $0 \leq j \leq \frac{m-3}{2}$ if $m$ is odd.  ({\it Cf.} 
\cite{Ba3}, Section 4.)  In particular, $z_{i-1,i}=1$ if $m \in 
\{2,3,4\}$.

We now use the fact that each vertex $s_i$ appears in two edges in $C$.  
Because $s'_i = z_{i-1,i}^{-1} w_i s_i w_i^{-1} z_{i-1,i}$ and $s''_i = 
w_{i+1} s_i w_{i+1}^{-1}$, both $s'_i$ and $s''_i$ are conjugate to $s_i$, 
and therefore to each other.  Let $P_i = \{ 
[s'_is_{i,1}],...,[s_{i,r}s''_i] \}$ be a path of minimal length from 
$s'_i$ to $s''_i$, all of whose edges have odd labels.  Then the element

$$\barv_i = v_{s_{i,1}s'_i} v_{s_{i,2}s_{i,1}} \cdots v_{s''_is_{i,r}}$$

\noindent conjugates $s''_i$ to $s'_i$.

We now compute:

$$w_i s_i w_i^{-1} = z_{i-1,i} s'_i z_{i-1,i}^{-1} = z_{i-1,i} \barv_i 
s''_i \barv_i^{-1} z_{i-1,i}^{-1} = z_{i-1,i} \barv_i w_{i+1} s_i 
w_{i+1}^{-1} \barv_i^{-1} z_{i-1,i}^{-1}.$$

\noindent Thus

$$w_i^{-1} z_{i-1,i} \barv_i w_{i+1} \in C(s_i) = C(w_i^{-1} z_{i-1,i} 
s'_i z_{i-1,i}^{-1} w_i) = w_i^{-1} z_{i-1,i} C(s'_i) z_{i-1,i}^{-1} 
w_i.$$

\noindent Finally, we obtain

$$w_{i+1}w_i^{-1} = \barv_i^{-1} \bars_i z_{i-1,i}^{-1} \eqno(1)$$

\noindent for some $\bars_i \in C(s'_i)$.  Denote by $x_i$ the word 
appearing on the right-hand side of (1).  Then

$$x_k x_{k-1} \cdots x_2 x_1 = w_1w_k^{-1} w_kw_{k-1}^{-1} \cdots 
w_2w_1^{-1} = 1.$$

If we choose a geodesic representation of $\bars_i$, each of the words 
$\barv_i^{-1}$, $\bars_i$, and $z_{i-1,i}^{-1}$ are as short as possible.  

Recalling the generators of $C(s'_i)$ given by Theorem~\ref{brink}, 
$\bars_i$ may terminate with a term of the form $u_{s''_{i-1}s'_i}$ (if 
$m=m_{s''_{i-1}s'_i}$ is even) or $v_{s'_is''_{i-1}}$ (if $m$ is odd).  In 
this case, we can reduce the product appearing in (1) by multiplying this 
term with $z_{i-1,i}^{-1}$ if $z_{i-1,i} \neq 1$.

There is one other place where reduction can occur in $x_i$.  Suppose 
$\bars_i$ begins with a word of the form

$$v_{\alpha_1s'_i} v_{\alpha_2\alpha_1} 
\cdots v_{\alpha_l\alpha_{l-1}} v_{s'_i\alpha_l}$$

\noindent corresponding to an odd ``loop'' based at $s'_i$ in the diagram 
$\sv'$.  If 

$$\barv_i = v_{s_{i,1}s'_i} v_{s_{i,2}s_{i,1}} \cdots v_{s''_is_{i,r}} 
\neq 1,$$

\noindent we may have

$$s_{i,1}=\alpha_l, s_{i,2}=\alpha_{l-1}, ... , s_{i,j}=\alpha_{l-j+1},$$

\noindent where $j \leq r$.  Therefore after cancellation $\barv_i^{-1} 
\bars_i$ begins with a word of the form

$$v_{\beta_1s''_i} v_{\beta_2\beta_1} \cdots v_{s'_i\beta_j} \eqno(2)$$

\noindent for some $j \geq 0$.  (When $j=0$, (2) has the form 
$v_{s'_is''_i}$.)

After all cancellation has been performed on $x_i$, we obtain a new 
product, of ``even'' words $E$ (involving terms $u_{ss'_i}$) and ``odd'' 
words $O$ (involving terms $v_{st}$):

$$x_i = v_{\alpha_1s''_i} v_{\alpha_2\alpha_1} \cdots v_{s'_i\alpha_j} E_1 
O_1 E_2 O_2 \cdots E_l O_l w(s'_i,s''_{i-1}) \eqno(3)$$

\noindent where $w(s'_i,s''_{i-1})$ is some word in the letters $s'_i$ and 
$s''_{i-1}$.  (We allow $E_1$ and $O_l$ to be trivial.)  If $s'_i$ 
appears in $\bar{s}_i$ as a generator of $C(s'_i)$, it may be absorbed 
by $w(s'_i,s''_{i-1})$.  The exact structure of the words $E$ and $O$ is 
governed by Theorem~\ref{brink}.

All of the terms $v_{st}$, and all of the even terms $u_{ss'_i}$ which 
consist of more than a single letter $s$ will be called {\it long terms}.  
In case $u_{ss'_i}$ is a long term of length $2r+1$ (for $r \geq 1$), we 
will also call the words $(ss'_i)^r$ and $(s'_is)^r$ long terms by a 
slight abuse of terminology.  Terms consisting of a single letter (either 
$s'_i$ or $s$ such that $ss'_i=s'_is$) will be called {\it short terms}.

We claim that the product given in (3) is in fact geodesic; this is 
shown by another application of the result of Tits.  Thanks to the form of 
the long terms and two-dimensionality of $W$, the only possible subwords 
in the right-hand side of (3) which admit replacement as in TP would come 
from $\bars_i$.  ({\it E.g.}, TP may allow us to bring two instances of 
$s'_i$ together, in order to cancel them.)  But we have assumed $\bars_i$ 
to be geodesic, and therefore unchanged under application of TP (such 
cancellations have already been performed).

We are now ready to begin the proof of Theorem~\ref{cmatch}, inducting 
upon the length $k$ of the circuit $C$.

\section{The base cases}

We continue to use the notation from the previous sections, and prove 
Theorem~\ref{cmatch} and Theorem~\ref{cmatch2} for cycles of lengths 3 and 
4.  Some of the methods used in this section will be generalized in the 
following section, and so will be stated in general terms.  We will use 
the fact that $x_k x_{k-1} \cdots x_1 = 1$ in order to show that each word 
$x_i$ must have a very specific form.  The form of $x_i$ will allow us 
both to identify a circuit in $\sv'$ to which $C$ corresponds as in 
Theorem~\ref{cmatch} and to prove the statements regarding $w_iw_j^{-1}$ 
made in Theorem~\ref{cmatch2}.

Hereafter we say that the words $x_i \neq 1$ and $x_j \neq 1$ ($i>j$) are 
{\it adjacent} in the product $x_k \cdots x_1$ if either $i=j+1$ or 
$x_l=1$ for $j<l<i$.

We first assume there is no cancellation of common short terms between two 
words $x_i$ and $x_j$.  Having proven the theorems in this case, we will 
then indicate how to prove the theorems in general.  (In fact, by 
two-dimensionality, the case in which $k=3$ yields very little such short 
term cancellation, as there can be no vertex $\alpha$ not on $C'$ such 
that $\alpha$ commutes with two distinct elements $s'_i$ and $s'_j$.)

Let us first consider the case of a circuit of length 3: $x_3x_2x_1=1$.  
Unless all three words $x_i$ are trivial (in which case $w_1=w_2=w_3$ 
already and $C$ clearly corresponds to a circuit $C'$ in $\sv'$), at least 
two of these words are nontrivial.

\vskip 3mm

\noindent {\bf Case 1.} \ Suppose that $x_1=1$ (after renaming, if 
necessary).  Thus $x_3x_2=1$ and $s'_1=s''_1$.

First suppose that $x_3$ ends with the long even term $\alpha 
(s'_3\alpha)^{\frac{m}{2}-1}$ and $x_2$ begins with the long even term 
$(\beta s'_2)^{\frac{n}{2}-1} \beta$.  (If $x_2$ were not to begin with an 
even term, at most one pair of letters would cancel, and an application 
of TP would yield a contradiction.)  This implies that $s'_2=s''_2$.  
Since $s'_3 \neq s''_2 = s'_2$, we can avoid the same contradiction only 
if $\alpha=s'_2$ and $\beta=s'_3$, in which case $m=n$ and 

$$\alpha (s'_3\alpha)^{\frac{m}{2}-1} \cdot (\beta s'_2)^{\frac{n}{2}-1} 
\beta = s'_3s'_2.$$

In this case, easy computations (and applications of TP) show that there 
is no further cancellation if $u_{s'_2s'_3}$ is either the first term in 
$x_3$ or is preceded by another long term, and if $u_{s'_3s'_2}$ is 
either the last term in $x_2$ or is followed by another long term.  
Therefore, in order that $x_3x_2=1$, $u_{s'_2s'_3}$ must be preceded by 
$s'_3$ in $x_3$, and $u_{s'_3s'_2}$ must be followed by $s'_2$ in $x_2$.  
Moreover, if there are any further terms in $x_3$ and $x_2$, there can be 
no further cancellation.  Thus $x_2=x_3=(s'_2s'_3)^{\frac{m}{2}}$, and 
$s'_i=s''_i$ for $i=1,2,3$.  Moreover, $x_1=1 \Rightarrow w_1=w_2$, and 
$(s'_2s'_3)^{\frac{m}{2}}$ commutes with both $s'_2$ and $s'_3$.  
Therefore, $w_1s_1w_1^{-1}=s'_1$, $w_1s_2w_1^{-1}=s'_2$, and 

$$w_1^{-1}s'_3w_1=w_3^{-1} (s'_2s'_3)^{\frac{m}{2}} s'_3 
(s'_2s'_3)^{\frac{m}{2}} w_3 = w_3^{-1}s'_3w_3 = s_3.$$

\noindent Thus the same element (namely, $w_1$) conjugates each $s_i$ to 
$\hat{s}_i=s'_i$, proving Theorem~\ref{cmatch2} for $C$.

Now suppose $x_3$ ends with a long odd term $(s'_3 
\alpha)^{\frac{m-1}{2}}$ and $x_2$ begins with a long odd term $(\beta 
s''_2)^{\frac{n-1}{2}}$.  (As before, we obtain a contradiction to 
$x_3x_2=1$ if one term is odd and the other even.)  In order that more 
than one pair of letters cancel, it must be that $\alpha=s''_2$ and 
$\beta=s'_3$, so that $m=n$ and

$$(s'_3 \alpha)^{\frac{m-1}{2}} \cdot (\beta s''_2)^{\frac{n-1}{2}} = 
s''_2s'_3.$$

As before there is no further cancellation possible if $v_{s'_3s''_2}$ is 
either the first term in $x_3$ or is preceded by a long term and 
$v_{s'_3s''_2}$ is either the last term in $x_2$ or is followed by a long 
term.  In fact, only if $s'_3=s'_2$ and $v_{s'_3s''_2}$ is followed in 
$x_2$ by the term $s'_2$ is there further cancellation.  However, we are 
still left with a stray letter $s''_2$, so this product cannot occur.

However, if instead $x_3$ ends with $v_{s'_3s''_2}s'_3$ and $x_2$ begins 
with $v_{s'_3s''_2}$, we obtain the product 
$s'_3v_{s''_2s'_3}v_{s'_3s''_2} = s'_3$ in the middle of $x_3x_2$.  No 
further cancellation is possible unless $s'_2=s'_3$ and $v_{s'_3s''_2}$ is 
followed in $x_2$ by the letter $s'_2$.  In this case, it is easily seen 
that there are no more terms in either $x_2$ or $x_3$, so that 
$s'_2=s'_3$, $s''_2=s''_3$, and $x_2=x_3=s'_3v_{s''_2s'_3}$.  As before, 
$w_1=w_2$, so $w_1s_1w_1^{-1}=s'_1$ and $w_1s_2w_1^{-1}=s'_2$.  Now, 
however, the vertex $\hat{s}_3$ to which $s_3$ is to be conjugated is not 
$s'_3$, but $s''_3 = v_{s'_3s''_3} s'_3 v_{s'_3s''_3}^{-1}$.  This is 
seen by drawing the circuit $C'$ (compare Lemma~\ref{whichvertex}).  But 
note 

$$w_1^{-1} s''_3 w_1 = w_3^{-1} v_{s'_3s''_3} s'_3 s''_3 s'_3 
v_{s'_3s''_3}^{-1} w_3 = w_3 s'_3 w_3^{-1} = s_3.$$

\noindent Therefore a single element again conjugates the vertices 
appropriately.  Schematically, the two possibilities above can be 
summarized respectively as $x_3=E$, $x_2=E^{-1}$; and $x_3=O$, 
$x_2=O^{-1}$ (as before, $E=$``even'', $O=$``odd'').

Similar arguments show that $z_{2,3}=z_{1,2}=1$, that $x_3$ cannot end 
with a short term, and that $x_2$ cannot begin with a short term.  Thus 
the two products shown above are the only valid possibilities when 
$x_1=1$.

\vskip 3mm

\noindent {\bf Case 2.} \ Now suppose that each $x_i \neq 1$, and that 
there is no cancellation of common short terms.  We give the possible 
forms for $x_i$ schematically (as was done above) in the proposition 
below, leaving precise computations to the reader.

\begin{proposition} \label{form3}
Let $x_i \neq 1$ for $i=1,2,3$.  Up to renumbering, one of the following 
holds:

\vskip 2mm

\begin{tabular}{ll}

$x_3 = E_1E_2$, $x_2 = E_2^{-1} E_3$, $x_1 = E_3^{-1}E_1^{-1}$; &
$x_3 = E_1E_2$, $x_2 = E_2^{-1}$, $x_1 = E_1^{-1}$; \\

$x_3 = E_1O_1$, $x_2 = O_1^{-1}O_2$, $x_1 = O_2^{-1}E_1^{-1}$; & 
$x_3 = E_1O_1$, $x_2 = O_1^{-1}$, $x_1 = E_1^{-1}$; \\

$x_3 = E_1O_1$, $x_2 = O_1^{-1}E_2$, $x_1 = E_2^{-1}E_1^{-1}$ &
$x_3 = O_1E_1$, $x_2 = E_1^{-1}$, $x_1 = O_1^{-1}$; \\

$x_3 = O_1O_2$, $x_2 = O_2^{-1}O_3$, $x_1 = O_3^{-1}O_1^{-1}$; &
$x_3 = O_1O_2$, $x_2 = O_2^{-1}$, $x_1 = O_1^{-1}$. \\

\end{tabular}

\vskip 2mm

\noindent Here, $E_j$ represents either $u_{\alpha\beta}$ or 
$u_{\alpha\beta} \beta$ and $O_j$ represents either $v_{\gamma\delta}$ or 
$v_{\gamma\delta} \gamma$, for the appropriate choices of $\alpha$, 
$\beta$, $\gamma$, and $\delta$.
\end{proposition}

Of course, one may not be able to choose freely whether $E_j$ represents 
$u_{\alpha,\beta}$ or $u_{\alpha,\beta} \beta$, and similarly for $O_j$.  
That is, the exact forms of the words $E_j$ and $O_j$ are clearly 
interdependent.  ({\it Cf.} Section 4 of \cite{Ba3}.)

Note that in multiplying any two distinct terms $x_i$ and $x_j$, at 
most one long term may cancel.  (That is, as we saw above, forms such as 
$x_1=1$, $x_3=O_1O_2$, $x_2=O_2^{-1}O_1^{-1}$ cannot occur.)  In fact, 
this will remain true even as we consider arbitrarily long circuits.

\vskip 2mm

\noindent {\bf Remark.} \ There are a few cases which must be handled 
carefully; these cases involve the affine Euclidean Coxeter groups 
whose diagrams are triangles with edge label multisets $\{2,3,6\}$, 
$\{2,4,4\}$, or $\{3,3,3\}$.  For instance, suppose that $x_i=s'_is''_i$ 
and $s'_i=s''_{i+1}$ for $i=1,2,3$.  Then $x_3x_2x_1=1$, 
$\hat{s}_i=s'_{i+1}$, and the edges of $\sv'$ which correspond to those in 
$C$ do indeed form a circuit of length 3.  However, it can be shown (with 
the aid of Lemma~\ref{edgenormalizers}) that no $w \in W$ satisfies 
$ws_iw^{-1}=\hat{s}_i$.  We claim that these forms of $x_i$ lead to a 
fundamental contradiction.  Consider the parabolic subgroup 
$w_1W_Cw_1^{-1}$; because $w_1s_1s_1^{-1}=s'_1$, $w_1s_3w_1^{-1}=s''_3$, 
and $w_1s_2w_1^{-1}=s'_3s'_1s''_3s'_1s'_3$, $w_1W_Cw_1^{-1} \subseteq 
W_{C'}$.  Now we apply the construction of $x_i$ given in Section 3 ``in 
reverse'', proceeding from $C'$ to $C$.  (Essentially, we compute the 
ratios $x'_i=w_{i-1}^{-1}w_i$.)  We have proven our result for $C$ (by 
conjugating back from the corresponding result for $C'$) provided the 
words $x'_i=w_{i-1}^{-1}w_i$ do not have forms similar to those of $x_i$.  
In particular, we may assume that $x'_i \in W_C$.  Therefore, since 
$s''_3=s'_2$ and $w_1^{-1}s'_3w_1=w_1^{-1}w_3s_3w_3^{-1}w_1 \in W_C$, 
$W_{C'} \subseteq w_1W_Cw_1^{-1}$.  Therefore $w_1W_Cw_1^{-1}=W_{C'}$.  
The set

$$\{w_1s_iw_1^{-1} \ | \ i=1,2,3\} = 
\{s'_1,s''_3,s'_3s'_1s''_3s'_1s'_3\}$$

\noindent must therefore generate $W_{C'}$.  However, it can be shown 
that this set does not generate $W_{C'}$ (this is possible because the 
affine Euclidean group $W_C$ is not cohopfian).  This gives a 
contradiction.

Therefore the case in which $x_i=s'_is''_i$ and $s''_i=s'_{i-1}$ for all 
$i$ cannot occur.  Any similar case involving the affine Euclidean 
triangle groups can be outlawed in an analogous fashion, and we are forced 
to conclude that the forms for $x_1 \cdot x_2 \cdot x_3$ given above are 
exhaustive.

\vskip 2mm

Before turning our attention to a proof of Theorem~\ref{cmatch2}, we state 
the following result concerning circuits of length 4.

\begin{proposition} \label{form4}
Let $C$ be a circuit of length 4 and let $x_4x_3x_2x_1=1$.  Assume 
furthermore that there is no cancellation of short terms between different 
words $x_i$ and $x_j$.  Then, up to a renumbering of the vertices, there 
are 27 forms for the product $x_4 \cdot x_3 \cdot x_2 \cdot x_1$ 
(entirely anologous to those given in Proposition~\ref{form3}).  In each 
case, every word $x_i$ has at most two long terms, every long term of 
$x_i$ cancels with a long term in either $x_{i+1}$ or $x_{i-1}$, and no 
more than two long terms cancel in any product $x_ix_{i-1}$.
\end{proposition}

One may check that in each of the cases mentioned, the exact forms of the 
words $O_i$ and $E_i$ which appear are determined completely, and that the 
subdiagram of $\sv'$ corresponding to these trivial products is a circuit 
of length 4 whose edges appear in the appropriate order.  This establishes 
Theorem~\ref{cmatch} for $k=4$.  (The form of the word which conjugates 
each edge $[s_is_{i+1}]$ appropriately is easy to compute, given the forms 
of $x_{i+1}$ and $x_i$.)

To prove Theorem~\ref{cmatch2} for $C$, we must first decide to which 
vertex $\hat{s}_i$ in $C'$ a given vertex $s_i$ in $C$ corresponds in 
Theorem~\ref{cmatch}.  To this end, we have the next lemma, which, along 
with the results which follow it, is stated in very general terms as it 
will be useful in the following section as well.  We say that two long 
terms {\it completely cancel} provided that they comprise the same two 
letters (thus their product has length at most 2).

\begin{lemma} \label{whichvertex}
Let $s_i$ be a vertex on $C$, and suppose that every word $x_l$ has one of 
the schematic forms $B$ (for ``blank''; i.e., no long terms), $O$, $E$, 
$OO$, $OE$, $EO$, or $EE$, for $1 \leq l \leq k$.  Suppose each long term 
completely cancels with a long term in an adjacent word, and no two words 
$x_i$ and $x_j$ allow complete cancellation of more than one long term 
between them.  Then $\hat{s}_i=s'_i$ unless $x_i$ has one of the forms 
$O$, $EO$, or $OO$ and the final long term in $x_i$ completely cancels 
with the first long term in $x_{i-1}$.  In this case, 
$\hat{s}_i=s''_{i-1}$.
\end{lemma}

\begin{proof}
The cases to consider depend on the schematic form of $x_i$.  If $x_i$ has 
no long terms, the lemma is clearly true.  We prove one of the nontrivial 
cases and leave the rest to the reader.

Suppose that $x_i$ has the form $OO$.  The second long term in $x_i$ must 
cancel with the first long term in $x_{i-1}$.  Therefore the long terms 
$v_{s'_i\alpha}$ and $v_{\beta s''_{i-1}}$ have the same letters, and 
$\alpha=s''_{i-1}$, $\beta=s'_i$.  Similarly, since the first long term 
in $x_i$ must completely cancel with the last in $x_{i+1}$, 
$\beta=s'_{i+1}$.  Therefore the edges $[s'_is''_{i-1}]$ and 
$[s'_{i+1}s''_i]$ meet at the vertex $s''_{i-1}=s''_{i+1}$.  But these are 
the edges of $C'$ which match the edges $[s_{i-1}s_i]$ and $[s_is_{i+1}]$ 
of $C$ as in Theorem~\ref{cmatch}.  Therefore the vertex $s''_{i-1}$ must 
be matched with $s_i$.
\end{proof}

We continue to consider arbitrary $k \geq 3$.  From the results of Section 
5 it will follow that provided there is no cancellation of short terms 
between different words $x_i$, every word $x_i$ will possess at most two 
long terms and have one of the seven schematic forms given above.  Also, 
if there is no cancellation of short terms, we claim that the circuit $C$ 
is strongly rigid in the sense that there is a single word $w$ which 
conjugates every vertex of $C$ appropriately.

To construct this word $w$, we must understand the form of $x_k$.  Given 
there is no short term cancellation, $x_k$ has either one or two long 
terms, perhaps followed by $s'_k$.  Suppose, for instance, that 
$x_k=v_{\alpha s''_k} v_{s'_k \alpha} s'_k$.  This forces $x_{k-1}$ to 
begin with the long term $v_{s'_ks''_{k-1}}$ (perhaps followed by 
$s'_k=s'_{k-1}$) and $x_1$ to end with the long term $v_{s'_1s''_k}$ 
(perhaps followed by $s'_1=\alpha$).  If indeed 
$x_{k-1}=v_{s'_ks''_{k-1}}$, then

$$w_1w_{k-1}^{-1}=x_kx_{k-1} = v_{\alpha s''_k} s'_k.$$

\noindent Since this product clearly does not conjugate $s_{k-1}$ to 
$\hat{s}_{k-1}=s'_{k-1}$ as needed, $w=w_1$ cannot be.  However, 
premultiplying by $\Delta_{s'_1s''_k}$ gives

$$\Delta_{s'_1s''_k} x_k x_{k-1} = 1,$$

\noindent and as one can check that $w=\Delta_{s'_1s''_k}w_1$ conjugates 
$s_1$ to $s''_k$ and $s_k$ to $s'_1=s''_{k-1}$, this $w$ works for the 
vertices $s_1$, $s_k$, and $s_{k-1}$.  We claim that such a $w$ can always 
be constructed in similar manner, depending on the precise form of $x_1$, 
$x_k$, and $x_{k-1}$, leaving the proof of this claim to the reader. 
({\it Cf.} Lemma~4.2 of \cite{Ba3}.)

\begin{lemma} \label{overthehump}
There exists a word $\pi \in 
\{1,s'_k,\Delta_{s'_1s''_k},\Delta_{s'_1s''_k}s'_k\}$ such that

$$\pi w_1 s_i w_1^{-1} \pi^{-1} = \hat{s}_i$$

\noindent for $i \in \{1,k-1,k\}$.
\end{lemma}

If $k=3$, this lemma completes the proof of Theorem~\ref{cmatch2}.  We 
claim that in fact this word $w$ will appropriately conjugate every vertex 
in $C$, even if $k \geq 4$.  ({\it Cf.} Proposition~4.3 of \cite{Ba3}.)

\begin{lemma} \label{whichword}
Let $k \geq 4$ and suppose that every word $x_i$ has one of the forms $B$, 
$O$, $E$, $OO$, $OE$, $EO$, or $EE$, for every $1 \leq i \leq k$.  
Further, suppose there is no short term cancellation between different 
words $x_i$.  Define $w=\pi w_1$ accordingly, as in 
Lemma~\ref{overthehump}.  Then $ws_iw^{-1}=\hat{s}_i$ for all $i$, where 
$\hat{s}_i$ is the vertex of $C'$ to which $s_i$ corresponds.
\end{lemma}

This lemma will complete the proof of Theorem~\ref{cmatch2} in case $k 
\geq 4$ and there is no cancellation of common short terms, assuming that 
every word $x_i$ has one of the seven schematic forms shown above.  
(Again, this last statement can be verified by computation in case $k=4$, 
and will be proven in Section 5 in case $k \geq 5$.)  In order to prove 
Lemma~\ref{whichword}, we make use of the following fact, which requires 
(when $k \geq 5$) an argument similar to, but simpler than, that used to 
prove Proposition~\ref{noiltc}.  (Essentially, 
Lemma~\ref{lotsocancellation} shows that once we have got over the initial 
``hill'' by premultiplying $w_1$ with $\pi$, there is sufficient 
cancellation between words $x_i$ to guarantee that the ratio $w_1w_i^{-1}$ 
is very short.)

\begin{lemma} \label{lotsocancellation}
Suppose that $k \geq 4$, that every word $x_i$ has one of the forms $B$, 
$O$, $E$, $OO$, $OE$, $EO$, or $EE$, and that there is no short term 
cancellation between different words $x_i$.  Define $t_i=\pi x_k \cdots 
x_i$ where $\pi$ is as above.  If the final long term of $x_i$ completely 
cancels with the first long term in $x_{i-1}$, then $t_i \in \{u,us'_i\}$, 
where $u$ is the final long term in $x_i$; otherwise $t_i \in \{1,s'_i\}$.  
\end{lemma}

\begin{proof}
We prove the case in which $x_k=v_{s'_ks''_k} s'_k$ and this term 
completely cancels with the last term in $x_1$ (the other cases will be 
similar).  In this case, $\pi=\Delta_{s'_1s''_k}$ and 
$t_i=\Delta_{s'_1s''_k}x_k \cdots x_i$ for all $i$.  The lemma is clearly 
true in case $i=k$.  Suppose that we have established the lemma for some 
fixed value of $i$, and consider $t_{i-1}$.

Suppose first that the last long term in $x_i$ is the odd term $v_{s'_i 
\alpha}$, and that this term completely cancels with the first long term 
in $x_{i-1}$.  First let $x_{i-1}=v_{s'_{i-1}s''_{i-1}} \epsilon_{i-1}$, 
$\epsilon_{i-1} \in \{1,s'_{i-1}\}$.  Then $t_i=v_{s'_i \alpha} s'_i 
\Rightarrow t_{i-1}=v_{s'_i \alpha} s'_i \cdot v_{s'_{i-1}s''_{i-1}} 
\epsilon_{i-1}$.  Complete cancellation implies $s'_i=s'_{i-1}$ and 
$\alpha=s''_{i-1}$, so $t_{i-1}=s'_{i-1} \epsilon_{i-1} \in 
\{1,s'_{i-1}\}$, as needed.  If instead 
$x_{i-1}=v_{s'_{i-1}s''_{i-1}}u_{\beta s'_{i-1}} \epsilon_{i-1}$, 
$t_{i-1}=s'_{i-1} u_{\beta s'_{i-1}} \epsilon_{i-1} \in \{ u_{\beta 
s'_{i-1}},u_{\beta s'_{i-1}}s'_{i-1}\}$, as needed.  Similar computations 
hold in case $x_{i-1}=v_{\beta s''_{i-1}} v_{s'_{i-1} \beta} 
\epsilon_{i-1}$.

If $t_i=v_{s'_i \alpha}$ and $x_{i-1}=v_{s'_{i-1}s''_{i-1}} 
\epsilon_{i-1}$, we compute $t_{i-1} \in \{s''_{i-1},s''_{i-1}s'_{i-1}\}$.  
It is at this point that a proof like that for Proposition~\ref{noiltc} 
must be used, in order to show that $s''_{i-1}$ cannot appear in $x_k 
\cdots x_{i-1}$.  (One must show that such an occurrence of $s''_{i-1}$ 
cannot be canceled by completing the product $t_{i-1} \cdot x_{i-2} \cdots 
x_1$.  This cancellation would require another occurrence of $s''_{i-1}$ 
in $x_{i-2} \cdots x_1$.  But if there were such an occurrence, we would 
be able to apply Theorem~\ref{cmatch2} inductively in order to obtain a 
contradiction, as in the proof of Proposition~\ref{noiltc}.)  Therefore 
this case cannot really occur.  The cases in which $t_i=v_{s'_i \alpha}$ 
and $x_{i-1}$ is either $v_{s'_{i-1}s''_{i-1}} u_{\beta s'_{i-1}} 
\epsilon_{i-1}$ or $v_{\beta s''_{i-1}} v_{s'_{i-1} \beta} \epsilon_{i-1}$ 
are outlawed by similar arguments.

We leave the proofs of the remaining cases to the reader.
\end{proof}

As a consequence, $t_is'_it_i^{-1}=\hat{s}_i$ for all $i$.  
Lemma~\ref{whichword} now follows almost immediately, since now

$$ww_i^{-1}=\pi w_1w_i^{-1} = \pi x_k \cdots x_i = t_i$$

\noindent for all $i$.

We now address the issue of short term cancellation.  Note that not only 
may such cancellation occur initially (before any long terms have been 
canceled), it might also occur after complete cancellation of two or more 
long terms has been performed.  For example, let $x_4=s'_1s'_4s'_1 
\alpha_1 \cdots \alpha_l s'_3s'_4s'_3$, $x_3=(s'_4s'_3)^2$, $x_2=\alpha_l 
\cdots \alpha_1$, and $x_1=(s'_4s'_1)^2$, where $s'_1s'_4$ and $s'_3s'_4$ 
have order 4 and for every $i$, $\alpha_is'_j=s'_j\alpha_i \Leftrightarrow 
j \in \{2,4\}$.  Then after canceling the last long term in $x_4$ with the 
first (and only) term in $x_3$, we may cancel every $\alpha_i$ in the 
product $x_4x_3 \cdot x_2$.

This poses no significant problems, as we have designed 
Theorem~\ref{cmatch2} to handle the possibility that such cancellation 
occurs.  Returning to the above example, define 
$w=\Delta_{s'_1s''_k}w_1=(s'_1s'_4)^2 w_1$.  Then $ws_1w^{-1}=s'_1$ 
and $ws_4w^{-1}=s'_4$, but $ws_3w^{-1}=\alpha_1 \cdots \alpha_l s'_3 
\alpha_l \cdots \alpha_1$.  However, every $\alpha_i$ satisfies the 
separation condition of Theorem~\ref{cmatch2}, relative to $s'_3$ and 
$s'_1$.

Let $k=4$.  If short term cancellation occurs at any stage in multiplying 
the words $x_i$ together, arguments like those used to handle the case 
$k=3$ may be used to prove Proposition~\ref{form4}, where short terms may 
now be inserted in between the schematic long terms.  The same arguments 
show that there is no cancellation between a short term in $x_i$ and a 
long term in $x_{i-1}$ without additional cancellation of a long term from 
$x_{i-1}$.  (This sort of outlawed cancellation will be the focus of much 
of Section 5.)  That is, short terms must ultimately cancel with other 
short terms, and long terms with long terms.

Therefore the only problem short term cancellation poses arises 
when completing the proof of Theorem~\ref{cmatch2}.  Yet all of the 
arguments from Lemma~\ref{whichvertex} through 
Lemma~\ref{lotsocancellation} still hold, with slight modification, in 
case there is short term cancellation.  First, we give the analogue of 
Lemma~\ref{lotsocancellation}; the proof of 
Lemma~\ref{lotsocancellation2} is very similar.

\begin{lemma} \label{lotsocancellation2}
Suppose that $k>3$ and that every word $x_i$ has at most 2 long terms (we 
still have the same seven schematic structures).  Define $t_i$ as in 
Lemma~\ref{lotsocancellation}.  Then $t_i$ can be written 
geodesically as $\alpha_1 \cdots \alpha_l u \epsilon_i$, where

\vskip 2mm

\noindent 1. \ $u$ is the last long term in $x_i$ if this long term 
completely cancels with the first long term in $x_{i-1}$ and $u=1$ 
otherwise,

\vskip 2mm

\noindent 2. \ $\epsilon \in \{1,s'_i\}$, and

\vskip 2mm

\noindent 3. \ for each letter $\alpha$ either $\alpha \in \{s'_1,s'_i\}$ 
or $\alpha \not\in C'$ and $\alpha$ commutes with some $s'_j$, $i \leq j 
\leq k$.  If $\alpha=s'_1$, $[s_2s_1]$ and $[s_1s_k]$ have label 2, and if 
$\alpha=s'_i$, $[s_{i-1}s_i]$ and $[s_is_{i+1}]$ have label 2.
\end{lemma}

Compare the letters $\alpha$ in this lemma with the letters $\alpha$ that 
arise in Theorem~\ref{cmatch2}.

As when proving Lemma~\ref{lotsocancellation}, a proof along the lines of 
Proposition~\ref{noiltc} is required in order to show that $s'_1$ and 
$s'_i$ are the only letters of $C'$ which can arise in as letters 
$\alpha$, and only under the circumstances indicated.  Consider the 
following example, which suggests why this should be true.

Let

\vskip 2mm

\begin{tabular}{lll}

$x_9 = (\alpha s'_1)^3 \alpha \alpha_1 \alpha_2 (\beta s'_9)^2$, & $x_8 = 
(\gamma s'_8)^2$, & $x_7 = \alpha_2$, \\

$x_6 = s'_6s''_6s'_6$, & $x_5 = s'_5s''_5s'_5 \alpha_1$, & $x_4 = (\delta 
s'_3)^5 \delta$, \\

$x_3 = 1$ & $x_2 = (\delta s'_3)^5 \delta (\alpha s'_1)^3$, & $x_1 = 1$, 
\\

\end{tabular}

\vskip 2mm

\noindent where $\alpha$, $\alpha_1$, $\alpha_2$, and $s'_1$ commute with 
$s'_9$, $\alpha_2$ commutes with $s'_7$, $\alpha_1$ commutes with 
$s'_5=s'_6$, both $\delta$ and $s'_3$ commute with both $s'_4$ and 
$s'_2$, and $\alpha$ and $s'_1$ both commute with $s'_2$.  If 
further $\beta=s'_8$, $\gamma=s'_9$, $s''_5=s''_6$, then $x_9 \cdots 
x_1 = 1$.  We compute $w_1w_8^{-1} = x_9x_8 = (\alpha s'_1)^3 \alpha 
\alpha_1 \alpha_2$, which demonstrates how $s'_i$ may be present in 
$w_iw_j^{-1}$.  Let us also examine $w_1w_3^{-1}$ and $w_1w_2^{-1}$.

$$w_1w_3^{-1}=x_9 \cdots x_3 = (\alpha s'_1)^3 \alpha (\delta s'_3)^5 
\delta,$$

\noindent giving us an example of a ratio $w_iw_j^{-1}$ in which $s'_j$ 
appears.  (Notice that since $\delta$ and $s'_3$ do {\it not} commute, 
this ratio is geodesic as written.)  However, in order that $x_k \cdots 
x_1=1$ hold, we must cancel the word $(\delta s'_3)^5 \delta$; the only 
other letter which can commute with both $\delta$ and $s'_3$ is $s'_2$, 
and we see that $w_1w_2^{-1} = x_9 \cdots x_2 = 1$, as needed.

Assume now that short terms cancel with short terms, long terms cancel 
with long terms, and that no word $x_i$ has more than 2 long terms.  
(Again, all of these statements follow from direct computation if $k=4$ 
and from the results of Section 5 if $k \geq 5$).  We can now prove 
Theorem~\ref{cmatch2} in case $k \geq 4$.

Reindex, as needed, so that $i=k$.  Defining $w$ as in 
Lemma~\ref{whichword}, $ws_1w^{-1}=\hat{s}_1$ and $ws_kw^{-1}=\hat{s}_k$ 
still hold.  As in the proof of Lemma~\ref{whichword}, $w=t_iw_i$ for 
every $i$.  Theorem~\ref{cmatch2} now follows from 
Lemma~\ref{lotsocancellation2}.

\section{The inductive step}

We have now indicated proofs of Theorem~\ref{cmatch} and 
Theorem~\ref{cmatch2} for achordal circuits of length at most 4.  
Inductively, assume that we have established these theorems for all 
achordal circuits of length less than or equal to $k-1$, and consider an 
achordal circuit $C$ of length $k$ in $\sv$.  For each edge $[s_{i-1}s_i]$ 
in $C$, Theorem~\ref{ematch} provides an edge $[s''_{i-1}s'_i]$ to which 
$[s_{i-1}s_i]$ corresponds.

Much as in the previous section, we will multiply the terms $x_i$ 
together, one at a time, performing all possible cancellation and length 
reduction as we go.  Also, as before, we begin by assuming that there is 
no cancellation of short terms between different words $x_i$.

Our first lemma can be proven using the results from \cite{BrHo}.

\begin{lemma} \label{edgenormalizers}
Let $(W,S)$ be a 2-d Coxeter system with diagram $\sv$, and let $[st]$ be 
an edge in $\sv$.  Let $w \in W$ satisfy $\{wsw^{-1},wtw^{-1}\} = 
\{s,t\}$.

\vskip 2mm

\noindent 1. \ If $[st]$ has label 2, then $w \in \{1,s,t,st\}$.

\vskip 2mm

\noindent 2. \ If $[st]$ has label greater than 2, then $w \in 
\{1,\Delta_{st}\}$, where $\Delta_{st}$ is the longest element in 
$W_{\{s,t\}}$.  In case the label on $[st]$ is odd, then $w=1 
\Leftrightarrow wsw^{-1}=s \Leftrightarrow wtw^{-1}=t$.
\end{lemma}

Now we prove a technical lemma that will often be used to reduce our 
problem to a case already considered.

\begin{lemma} \label{starpiecing}
Let Theorem~\ref{cmatch} and Theorem~\ref{cmatch2} both be proven for 
achordal circuits of length at most $k-1$, and let $C$ be an achordal 
circuit in $\sv$ of length $k$.  Suppose that there is a vertex $\alpha$ 
adjacent to the vertices $s_{i_1},s_{i_2},...,s_{i_r}$ so that each of the 
circuits $\{[\alpha s_{i_l}],[s_{i_l}s_{i_l+1}],...[s_{i_{l+1}}\alpha]\}$ 
is achordal and of length less than $k$, for $1 \leq l \leq r$.  Then 
Theorem~\ref{cmatch} and Theorem~\ref{cmatch2} both hold for $C$ as well.
\end{lemma}

\begin{proof}
Reindexing, we let $l=i_l$ for $1 \leq l \leq r$, and let $C_l$ denote the 
circuit $\{[\alpha s_l],...,[s_{l+1}\alpha]\}$.  By hypothesis, to each 
$C_l$ there is a circuit $C'_l$ in $\sv'$ which corresponds, edge by edge, 
to $C_l$.  (As usual, we use ``prime'' notation to indicate the 
corresponding vertices.)  Moreover, the ratios $w_iw_j^{-1}$ of the 
elements which conjugate the edges of a given $C_l$ are governed by 
Theorem~\ref{cmatch2}.

Consider the edge $[s_{l+1}\alpha]$, lying in the circuits $C_l$ and 
$C_{l+1}$.  Applying Theorem~\ref{cmatch} to these circuits yields group 
elements $\hat{w}_l$ and $w_{l+1}$ such that

$$\{ \hat{w}_ls_{l+1}\hat{w}_l^{-1} , \hat{w}_l\alpha \hat{w}_l^{-1}\} = 
\{ w_{l+1}s_{l+1}w_{l+1}^{-1} , w_{l+1}\alpha w_{l+1}^{-1} \}.$$

\noindent Lemma~\ref{edgenormalizers} implies that $w_{l+1}^{-1}\hat{w}_l 
\in \{1,s_{l+1},\alpha,s_{l+1}\alpha\}$ if $[s_{l+1}\alpha]$ has label 2 
and $\hat{w}_l^{-1}w_{l+1} \in \{1,\Delta_{s_{l+1}\alpha}\}$ if 
$[s_{l+1}\alpha]$ has label greater than 2.  Since $\hat{w}_l \alpha 
\hat{w}_l^{-1} = \alpha'$ and $\hat{w}_l s_{l+1} \hat{w}_l^{-1} = 
s'_{l+1}$, we obtain $\hat{w}_lw_{l+1}^{-1} \in 
\{1,s'_{l+1},\alpha',s'_{l+1}\alpha'\}$ if $[s_{l+1},\alpha]$ has label 2 
and $\hat{w}_lw_{l+1}^{-1} \in \{1,\Delta_{s'_{l+1}\alpha'}\}$ if 
$[s_{l+1}\alpha]$ has label greater than 2.

Now, from Theorem~\ref{cmatch2} applied to $C_l$, $w_l\hat{w}_l^{-1}$ can 
be written as a product $\beta_1 \cdots \beta_p$, where for every letter 
$\beta_i$, either $\beta_i=s'_l$ (in which case $[s_l\alpha]$ has 
label 2), $\beta_i=s'_{l+1}$ (in which case $[s_{l+1}\alpha]$ has label 
2), or $\beta_i$ does not lie on $C'$.  (In this last case, 
Theorem~\ref{cmatch2} shows that $\beta_i$ does not lie on $C'_l$; 
if $\beta_i$ were to lie on $C'_{l'}$, $l' \neq l$, we would contradict 
the achordality of $C$.  For later use, we note that by the separation 
condition of Theorem~\ref{cmatch2}, each such $\beta_i$ must in fact 
commute with $\alpha$.)

We can now compute $w_lw_{l+1}^{-1} = 
w_l\hat{w}_l^{-1}\hat{w}_lw_{l+1}^{-1}$ for each $l$, $1 \leq l \leq r$.  
The product of these words, taken in order, is trivial in $W$.  Therefore, 
application of the Tits Process (TP) must yield the trivial word.

Assume that for some edge $[s_{l+1}\alpha]$ with label greater than 2, 
$\hat{w}_lw_{l+1}^{-1}=\Delta_{s'_{l+1}\alpha'}$.  One easily sees that 
the letter $s'_{l+1}$ appears only in this word when forming the product 
$w_1w_2^{-1} \cdots w_rw_1^{-1}$ as above, and then only in 
$\Delta_{s'_{l+1}\alpha'}$.  (In particular, $s'_{l+1}$ cannot arise as 
a letter $\beta_i$ in $w_l\hat{w}_l^{-1}$, because we have assumed 
$[s_{l+1}\alpha]$ has label exceeding 2.)  In applying TP, no two 
occurrences of the letter $s'_{l+1}$ are brought next to one another.  
This letter can therefore not be canceled, contradicting the product's 
triviality.  Therefore, $\Delta_{s'_{l+1}\alpha'}$ cannot appear.  In 
particular, $\hat{w}_ls_{l+1}\hat{w}_l^{-1} = w_{l+1}s_{l+1}w_{l+1}^{-1}$ 
and $\hat{w}_l\alpha \hat{w}_l^{-1} = w_{l+1}\alpha w_{l+1}^{-1}$ hold, so 
that there is no ``twisting'' at the edge $[s'_{l+1}\alpha']$ in $\sv'$ 
when the two adjacent circuits $C'_l$ and $C'_{l+1}$ are met along this 
edge.  This concludes the proof of Theorem~\ref{cmatch} for $C$.

Now for Theorem~\ref{cmatch2}.  Note first that the argument from the 
previous paragraph also shows that $\hat{w}_lw_{l+1}^{-1} \in 
\{1,\alpha'\}$ if the label on $[s_{l+1}\alpha]$ is 2.

Now take any two distinct vertices on $C$ and consider the edges 
$e_i=[s_is_{i+1}]$ and $e_j=[s_{j-1}s_j]$ in $C$ (as in 
Theorem~\ref{cmatch2}).  Suppose $e_i$ and $e_j$ lie on circuits $C_l$ 
and $C_m$, respectively.  We suppose $j>i$, and let $\bar{w}_i$, 
$\bar{w}_j$ be the elements conjugating $e_i$ and $e_j$ (with respect to 
$C_l$ and $C_m$, resp.) provided by Theorem~\ref{cmatch2}.  Then

$$\bar{w}_i\bar{w}_j^{-1} = \bar{w}_i\hat{w}_l^{-1} \cdot 
\hat{w}_lw_{l+1}^{-1} \cdots \hat{w}_{m-1}w_m^{-1} \cdot 
w_m\bar{w}_j^{-1}.  \eqno(4)$$

From the preceding arguments we know this product can be written in the 
letters $\beta$ which arise from applying Theorem~\ref{cmatch2} to each 
achordal cycle $C_p$ in turn, as well as the letter $\alpha'$.  We will 
have proven Theorem~\ref{cmatch2} for $C$ once we show each letter $\beta$ 
is either in $\{s'_i,s'_j\}$ or commutes with at least two letters in $C'$ 
and satisfies the separation condition of Theorem~\ref{cmatch2}.

First suppose that $\beta$ lies on $C'$ and neither $\beta=s'_i$ nor 
$\beta=s'_j$ holds.  Then $\beta=s'_p$ for some $p$ and $[s'_p\alpha]$ is 
an edge labeled 2.  The only two words in (4) which may contain the letter 
$s'_p$ are $\hat{w}_pw_{p+1}^{-1}$ and $w_{p+1}\hat{w}_{p+1}^{-1}$.  
Therefore, in applying TP to the product

$$\bar{w}_i\bar{w}_j^{-1} \cdot \bar{w}_j\hat{w}_m^{-1} \cdots 
\hat{w}_{l-1}w_l^{-1} \cdot w_l\bar{w}_i^{-1} = 1  \eqno(5)$$

\noindent the instances of $s'_p$ in $\hat{w}_pw_{p+1}^{-1}$ cancel with 
those in $w_{p+1}\hat{w}_{p+1}^{-1}$.  We claim that the same cancellation 
takes place in (4) as in (5) (that is, we do not need the context of (5) 
in order to use TP to cancel these occurrences of $s'_p$).  This is 
so because none of the letters that appear in between occurrences of 
$s'_p$ are related to one another, except $\alpha$.  (This follows from 
two-dimensionality and the fact that every $\beta \neq \alpha$ commutes 
with $\alpha$ when it appears in one of the two terms given above.)  
Therefore, we can perform this cancellation of letters $s'_p$ in (4), and 
$s'_p$ does not in fact occur.

Now suppose $\beta \not\in C'$.  The separating condition is satisfied 
unless $\beta$ is adjacent by edges labeled 2 only to $\alpha'$ and to a 
single letter $s'_p$ in some cycle $C'_{l'}$, $l \leq l' \leq m$.

Yet we claim that $\beta$ must also be adjacent to some other letter 
$s'_q$ in $C'$.  To see this, examine the product in (5) again.  As this 
word represents the trivial element in $W$, applying TP must yield 
the empty word.  If a letter $\beta$ occurs in $\bar{w}_i\bar{w}_j^{-1}$, 
it either cancels with another occurrence of $\beta$ in that same 
subproduct (in which case it can be removed from that subproduct 
altogether anyway, by an application of TP to the subproduct, as indicated 
above) or it cancels with an occurrence of $\beta$ lying in 
$\bar{w}_j^{-1} \cdots \bar{w}_i$.  This other occurrence of $\beta$ comes 
from an application of Theorem~\ref{cmatch2} to yet another cycle 
$C_{l''}$ which lies between $C_l$ and $C_m$, ``opposite'' the cycles 
$C_{l+1},...,C_{m-1}$.  Therefore there is a vertex $s'_q$ such that 
$[s'_q\beta]$ is an edge labeled 2 and such that the path 
$\{[s'_p\beta],[\beta s'_q]\}$ divides $C'$ into two cycles, one 
containing $s'_i$, and the other containing $s'_j$, as desired.

\end{proof}

The following proposition shows that there cannot be a great deal of 
cancellation between words $x_i$ and $x_j$ for which $d(i,j)$ is large. 
The arguments in this proof demonstrate the flavor of many of the 
arguments to come.

\begin{proposition} \label{xicancel}
Let Theorem~\ref{cmatch} and Theorem~\ref{cmatch2} both be proven for 
achordal circuits of length at most $k-1$, and let $C$ be an achordal 
circuit in $\sv$ of length $k$.  Define the words $x_i$ as above.  Let 
$x_i$, $x_j$, and $x_l$ be adjacent words in $x_k \cdots x_1$ (for 
$i>j>l$), and assume that there is no short term cancellation.

\vskip 2mm

\noindent 1. \ If $d(i,j) \geq 3$, then we may assume that $x_ix_j$ is a 
geodesic word for the group element it represents (i.e., no length 
reduction is possible).

\vskip 2mm

\noindent 2. \ If $d(i,j)=2$, then at most one pair of letters cancels in 
forming the product $x_ix_j$.  

\vskip 2mm

\noindent 3. \ Suppose $d(i,j)=2$ and $d(j,l)=2$.  If $x_j$ consists 
either of the single long term $v_{s'_js''_j}$ or the single long term 
$z_{j-1,j}^{-1} \neq 1$, then at least one of the products $x_ix_j$ or 
$x_jx_l$ allows no cancellation.
\end{proposition}

Note that if $i=k$ in the first case, then $d(i,j) \geq 3$ forces $j \neq 
1$, $j \neq 2$, and $i=k-1$ forces $j \neq 1$.  The third case will allow 
an easy proof of Theorem~\ref{cmatch} and Theorem~\ref{cmatch2} in case 
all adjacent $x_i$ and $x_j$ words satsify $d(i,j) \geq 2$.

\begin{proof}
1. \ First suppose that $x_i$ and $x_j$ are adjacent and $d(i,j) \geq 3$.  
Assume (to derive a contradiction) that there is length reduction in the 
product $x_ix_j$.  There are a few possibilities.

First, $x_i$ could end with a long term $(s'_i \alpha)^l$ while $x_j$ 
begins with a long term $(\beta s''_j)^{r}$ (of course, $s'_j=s''_j$ 
if $\barv_j=1$).  That there is reduction implies either that $\beta=s'_i$ 
and $\alpha=s''_j$ both hold or that $\alpha=\beta$.

In the first case, $s'_i$ and $s''_j$ are adjacent in $\sv'$.  Because 
$x_p=1$ for $j<p<i$, $\{[s'_is'_{i-1}],...,[s'_{j+1}s''_j]\}$ forms a 
circuit $D'$ of length at most $k-1$ in $\sv'$.  If $D'$ is not achordal, 
we can find a smaller circuit containing a subset of these vertices which 
is achordal in $\sv'$.  By inductive hypothesis, there is a circuit $D$ in 
$\sv$ which corresponds to $D'$ as in Theorem~\ref{cmatch}.  However, all 
but one of the edges of $D'$ correspond as in Theorem~\ref{ematch} to 
edges of $C$; therefore $D$ is a circuit whose vertices are vertices of 
$C$.  This implies that $C$ was not achordal, a contradiction.

In the second case, there is a vertex, $\alpha$, adjacent to both $s'_i$ 
and to $s''_j$.  As in the previous paragraph, we find a circuit $D'$
$\{[s'_is'_{i-1}],...,[s'_{j+1}s''_j],[s''_j \alpha],[\alpha s'_i]\}$ of 
length at most $k-1$.  If $D'$ is not achordal, it can be subdivided into 
shorter achordal circuits by addition of edges $[\alpha s'_p]$, $j<p<i$.  
By inductive hypothesis, each of these achordal circuits corresponds as 
in Theorem~\ref{cmatch} to an achordal circuit of the same length in 
$\sv$.  In fact, we can ``piece together'' these individual circuits by 
pasting along the common chords which correspond as in 
Theorem~\ref{ematch} to the edges $[\alpha s'_p]$.  That is, if $D'_1$ 
and $D'_2$ are two achordal circuits into which $D'$ has been divided and 
which share the edge $[\alpha s'_p]$, then the corresponding circuits 
$D_1$ and $D_2$ in $\sv$ share an edge $[\gamma s_p]$ corresponding to 
$[\alpha s'_p]$.  Moreover, because we know that the edges $[s_qs_{q-1}]$ 
in $\sv$ corresponding to $[s'_qs'_{q-1}]$ follow each other in sequence, 
there can be no ``twisting'' at the common edge $[\gamma s_p]$.

Now consider the circuit 
$\{[s_js_{j-1}],...,[s_{i+1}s_i],[s_i\gamma],[\gamma s_j]\}$ in $\sv$.  If 
this circuit is not achordal, by the achordality of $C$ it can be 
subdivided into shorter achordal circuits by addition of edges $[\gamma 
s_p]$.  Therefore $C$ itself has been subdivided into shorter achordal 
circuits, yielding the configuration described in Lemma~\ref{starpiecing}.  
We apply this lemma and conclude that Theorem~\ref{cmatch} and 
Theorem~\ref{cmatch2} hold for $C$, as desired.

The other possibilities for cancellation between $x_i$ and $x_j$ (i.e., 
a short term of $x_i$ canceling with a long term of $x_j$, or {\it vice 
versa}, or cancellation between long terms of different sorts, arising 
when $z_{i-1,i} \neq 1$, for instance) can be handled in an entirely 
similar fashion.  

\vskip 2mm

\noindent 2. \ Now assume that $x_i$ and $x_j$ are adjacent and that 
$d(i,j)=2$.  As before, let us first consider the case of a long term 
$(s'_i\alpha)^l$ in $x_i$ canceling with a long term $(\beta s''j)^r$ of 
$x_j$.  If there is to be more than one pair of letters canceling in 
$x_ix_j$, it must be that $\alpha=s''_j$ and $\beta=s'_i$.  As before, 
this forces $s'_i$ and $s''_j$ to be adjacent, yielding a shorter achordal 
circuit to which the inductive hypothesis can be applied, giving a 
contradiction to the achordality of $C$.

As in the first case, we leave the similar arguments for the other 
possibilities to the reader.

\vskip 2mm

\noindent 3. \ Suppose that $x_j=v_{s'_js''_j}$, that $x_i$ ends with 
the long term $(s'_i\alpha)^m$, and that $x_l$ begins with the long term 
$(\beta s''_l)^r$.  (The other cases for $x_j=v_{s'_js''_j}$ and the cases 
in which $x_j=z_{j-1,j}^{-1} \neq 1$ will be left to the reader.)

Since $s'_i=s'_j$ and $s''_j=s''_l$ both lead to a contradiction to $C$'s 
achordality, cancellation can only occur if $\alpha=s'_j$ or 
$\beta=s''_j$.  Suppose that both of these equalities hold.  Then there 
are circuits 
$C'_1=\{[s'_i,s'_{i-1}],[s'_{i-1}s''_j],[s''_js'_j],[s'_js'_i]\}$ 
and $C'_2=\{[s''_js'_j],[s'_js'_{j-1}],[s'_{j-1}s''_l],[s''_ls''_j]\}$ in 
$\sv'$, sharing the edge $[s''_js'_j]$.  First assume that these circuits 
are achordal.  Then by Theorem~\ref{cmatch} there are corresponding 
circuits $C_1$ and $C_2$ in $\sv$.  Each of $C_1$ and $C_2$ contains two 
consecutive edges from $C$, and $C_1$ and $C_2$ share a common edge.  
Moreover, because we know the sequence of edges in $C$, we conclude that 
the circuit formed by replacing the edges $[s_is_{i-1}],...,[s_{l+1}s_l]$ 
by $[s_i\gamma]$ and $[\gamma s_l]$ is shorter than $C$; if it is not 
achordal, then it can be subdivided by addition of chords $[\gamma s_n]$ 
for some letters $s_n$, and an application of Lemma~\ref{starpiecing} 
brings our proof to a close (by demonstrating that the ``twist'' 
apparent at the edge $[s''_js'_j]$ could not in fact have been, giving a 
contradiction).

If $C'_1$ and $C'_2$ had not been achordal, we could have shortened them 
by introducing the requisite edge and applying the arguments of the 
previous paragraphs in order to reach the same conclusion.

\end{proof}

Suppose now that for any two adjacent words $x_i$ and $x_j$, $d(i,j) \geq 
2$.  (Here we may assume that there could be short term cancellation, but 
we suppose that all common short terms have been canceled.)  Suppose also
that for some $j$, $\barv_j \neq 1$.

If $x_j$ contains at least two long terms, then $x_j$ has length at least 
4.  Consider the product $x_ix_jx_l$ where $x_i$ and $x_l$ are the 
adjacent words nearest to $x_j$ on either side.  At most two letters of 
$x_j$ are canceled in this product.  Can there be more cancellation upon 
multiplying more words?  Suppose that $x_l$ consists of a single letter, 
$\beta$, which cancels in the product $x_jx_l$.  Then either $\beta=s'_l$ 
or $\beta s'_l = s'_l \beta$; in any case, $s'_l=s''_l$.  If $x_m$ is the 
next nontrivial word after $x_l$, the arguments used to prove part 1 of 
Proposition~\ref{xicancel} can be applied to show that there is no 
cancellation in multiplying $x_j\beta \cdot x_m$.  The same can be said 
for the product $x_ix_j$.

Ultimately we conclude that in the product $x_k \cdots x_1$, there are 
letters of $x_j$ which remain uncanceled after all reduction is performed 
(including application of TP).  Thus $x_j$ could not have had more than 
one term.

If instead $\barv_j \neq 1$ and this is the only long term appearing in 
$x_j$, part 3 of Proposition~\ref{xicancel} allows us to reach the same 
conclusion: some of $x_j$ must remain in multiplying $x_k \cdots x_1$.

Therefore, $\barv_i=1$, so $s'_i=s''_i$ for all $i$, $1 \leq i \leq k$, 
and the circuit $C'$ corresponds as in Theorem~\ref{cmatch}.

Does this matching satisfy Theorem~\ref{cmatch2} as well?  Note that an 
argument similar to that in the above paragraphs shows that $z_{i-1,i}=1$ 
must hold for all $i$, $1 \leq i \leq k$ and no word $u_{\alpha 
s'_i}$ can occur in any $x_i$.  Thus, the ratio 
$w_{i+1}w_i^{-1}$ can be written $\alpha_1 \cdots \alpha_{l_i} 
\epsilon_i$, where $\alpha_l s'_i = s'_i \alpha_l$ for $1 \leq l \leq l_i$ 
and $\epsilon_i \in \{1,s'_i\}$.  Suppose that $s'_i$ were to appear.  As 
$x_k \cdots x_1=1$, some other word, $x_j$, must contain $s'_i$, and 
therefore $s'_is'_j=s'_js'_i$.  This contradicts the achordality of $C'$.  
Therefore $\epsilon_i=1$ for all $i$.  At last, an argument similar to 
that used in establishing the ``separation condition'' of 
Theorem~\ref{cmatch2} in Lemma~\ref{starpiecing} shows that the same 
condition is satisfied in this case.

\vskip 2mm

We must now turn our attention to the case in which for some adjacent 
words $x_i$ and $x_j$, $d(i,j)=1$.  (We frequently assume that 
$d(i,j)=1$ for all adjacent words $x_i$ and $x_j$.  Whenever this is 
not so, our arguments are often made simpler, as the reader is invited to 
verify.)  The arguments we offer, though at times technical, are entirely 
analogous to those that have come before.

Until further notice we will assume that $z_{i-1,i}=1$ for all $i$ (our 
arguments below will show that this must be the case anyway) and that 
there is no short term cancellation, as in the previous section.  We 
distinguish between two putative types of cancellation between two 
adjacent words $x_i$ and $x_{i-1}$: {\it complete long term cancellation} 
(CLTC) and {\it incomplete long term cancellation} (ILTC).

CLTC occurs when $x_i$ ends with a long term (followed perhaps by $s'_i$) 
and $x_{i-1}$ begins with a long term in the same letters.  In this case 
the resulting product has at most two letters and lies in 
$\{1,s'_i,s''_{i-1},s''_{i-1}s'_i\}$ if the common long term is odd, 
and $\{1,s'_i,s'_{i-1},s'_is'_{i-1}\}$ if the common long term is even.

ILTC is the cancellation that can conceivably occur when a long term in 
$x_i$ (resp. $x_{i-1}$) cancels with either a short term, a single letter 
of a long term, or both, in $x_{i-1}$ (resp. $x_i$).  For instance, if 
$x_i$ ends with $v_{s'_i\alpha}$ and $x_{i-1}$ begins with $\alpha 
u_{\beta s'_{i-1}}$ (where $\alpha \neq s'_{i-1}$), the letters $\alpha$ 
cancel.  Further, it is possible that $\beta=s'_i$, so that one more pair 
of letters cancels.  There is no more cancellation than this, leaving at 
least two letters of $u_{\beta s'_{i-1}}$ intact, justifying the use of 
the term ``incomplete''.  There will be a small number of exceptional 
cases of ILTC in which an entire long term is canceled; these cases will 
more closely resemble CLTC in some respects.

We shall describe all possible forms of CLTC and ILTC shortly.

The rough course of our remaining argument is as follows.  For each $x_i$ 
there is a subdiagram of $\sv'$ (which we will call a {\it piece}) whose 
form can be derived from $x_i$.  In case $x_i$ has at most two long terms, 
the corresponding piece has one of the seven general forms examined in 
Section 4.  When two words $x_i$ and $x_{i-1}$ exhibit any sort of 
cancellation, we are given information about how to put the corresponding 
pieces together, and given a chain of consecutive ``short'' words (in a 
sense to be introduced below) $x_i,x_{i-1},...,x_j$, each of which cancels 
in some way with the previous one, a subdiagram of $\sv'$ emerges which 
resembles a segment in the circuit $C'$ whose existence we wish to 
establish.

Putting together pieces, we first show that we may assume ILTC does not 
occur.  We then indicate how the same argument may be used to show that 
$z_{i-1,i}=1$ for all $i$, and that no word $x_i$ comprises more than 2 
long terms.  At this point we will be in a position to prove 
Theorem~\ref{cmatch} for $C$ much as was done for short circuits in 
Section 4.  Then we will appeal to the final results from that section 
to complete the proof of Theorem~\ref{cmatch2} for $C$.

We call a word $x_i$ {\it terse} if it has no more than 2 long terms.  As 
we have seen, there are seven general forms for such words: $B$, $O$, $E$, 
$OO$, $OE$, $EO$, and $EE$, where short terms may be inserted in 
appropriate places.  Note that in cases $B$, $E$, and $EE$, $s'_i=s''_i$.  
The pieces which correspond to each of these terse word forms are depicted 
in Figure~1.

\begin{figure}[hbt]
	\begin{center}
		\input{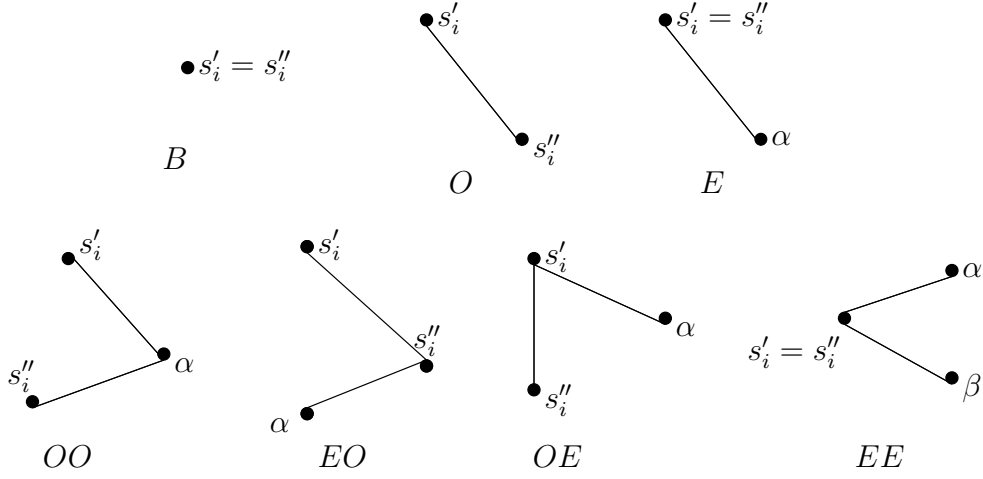}
	\end{center}
	\caption{Pieces for each of the seven forms of terse words}
\end{figure}

\begin{proposition} \label{noiltc}
We can assume that there is no incomplete long term cancellation between 
two words $x_i$ and $x_{i-1}$.
\end{proposition}

To begin our proof, let us assume that there is ILTC between $x_i$ and 
$x_{i-1}$; by reindexing, assume that $i=k$.

There may also be ILTC between $x_{k-1}$ and $x_{k-2}$, and then 
too between $x_{k-2}$ and $x_{k-3}$, and so forth.  We will continue to 
multiply terms $x_k,x_{k-1},...,x_j$ until we no longer see ILTC, keeping 
track of the product $x_kx_{k-1} \cdots x_j$ (and the corresponding 
concatenation of pieces) as we go.  Unless some nontrivial word remains 
when all words $x_i$ are multiplied, we will either be able to apply 
Theorem~\ref{cmatch} to a shorter circuit to obtain a contradiction, or be 
able to apply Lemma~\ref{starpiecing} to yield the desired conclusion by 
appeal to the inductive hypothesis.

For notational convenience, we denote the subdiagram of $\sv'$ formed by 
concatenating the pieces corresponding to $x_{i_1},x_{i_1-1},...,x_{i_2}$ 
by $\sv'(i_1,i_2)$.

We first multiply $x_k,x_{k-1}, \cdots, x_j$ only so long as each word 
$x_i$ for $j \leq i \leq k-1$ has at most one long term and there is no 
exceptional ILTC as defined below.  (For the time being the exact form of 
$x_k$ is not important.)  Thus at each step we multiply words of types 
$B$, $E$, and $O$ with one another.  Since we have assumed that there is 
no short term cancellation, and because no more than one short term 
may cancel with a long term (otherwise we would contradict the fact that 
$W$ is 2-d), we are left with the following possible products for $x_i 
\cdot x_{i-1}$:

\vskip 2mm

\begin{tabular}{lllll}

1. & $B \cdot O$, & $B \cdot E$, & $O \cdot B$, & $E \cdot B$, \\

2. & $O \cdot O$, & $O \cdot E$, & $E \cdot O$, & $E \cdot E$, \\

3. & $Os'_i \cdot O$, & $Os'_i \cdot E$, & $Es'_i \cdot O$, & $Es'_i \cdot 
E$, \\

4. & $O\alpha \cdot O$, \ (*) & $O\alpha \cdot E$, & $E\alpha \cdot O$, \ 
(*) & $E\alpha \cdot E$, \\

5. & $O \cdot \beta E$, \ (*) & $E \cdot \beta E$, & & \\

6. & $O \cdot s'_{i-1}E$, & $E \cdot s'_{i-1}E$. & & \\

\end{tabular}

\vskip 2mm

\noindent Here, $\alpha s'_i = s'_i \alpha$ and $\beta s'_{i-1} = s'_{i-1} 
\beta$.

The cases marked (*) represent cases in which the ILTC may be more 
complicated.  In each of these cases, one of the words $x_i$, $x_{i-1}$ 
can indeed be completely canceled in forming the product $x_ix_{i-1}$.  
For instance, if $x_i=v_{s'_i \gamma} \alpha$ (where $\alpha s'_i=s'_i 
\alpha$) and $x_{i-1}=v_{s'_{i-1}s''_{i-1}}$, then $\alpha=s'_{i-1}$ and 
$\gamma=s''_{i-i}$ can both hold.  If $s'_{i-1}s''_{i-1}$ has order $n=3$, 
then $x_{i-1}$ has been completely canceled.  (Incidentally, when $n=3$, 
two-dimensionality implies that $s'_i\gamma$ has order at least 7.  We use 
this fact later.)  However, one may apply Lemma~\ref{edgenormalizers} to 
show that further cancellation of the remainder of $x_i$ with $x_{i-2}$ 
would contradict the structure of $C$ (we would essentially be trying to 
``twist'' the diagram at an edge labeled 2, which cannot be done, by 
Lemma~\ref{edgenormalizers}).  Therefore we can cancel no more.  Because 
$x_i$ and $x_{i-2}$ do not cancel with each other, the chain of ILTC ends 
at this point, and may or not begin ILTC anew with $x_{i-2} \cdot 
x_{i-3}$.  If a long term {\it is} completely canceled in the product $x_i 
\cdot x_{i-1}$, we call this case of ILTC {\it exceptional}.

As an illustration of the above products, Figure~2 displays the 
subdiagrams $\sv'(i,i-1)$ corresponding to the products $O \cdot O$, 
$Os'_i \cdot E = v_{s'_is''_i}s'_i \cdot u_{\alpha s'_{i-1}}$, and $E 
\cdot \beta E = u_{\alpha s'_i} \cdot \beta u_{\gamma s'_{i-1}}$.

\begin{figure}[hbt]
	\begin{center}
		\input{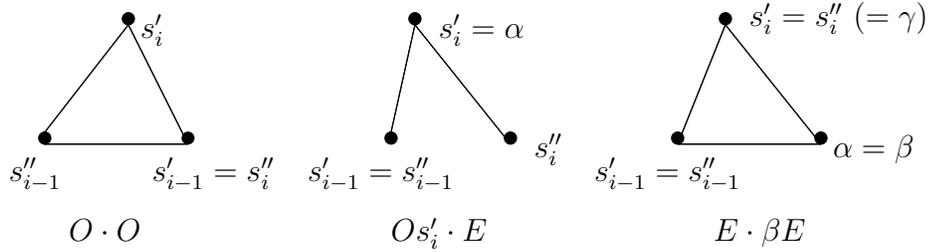}
	\end{center}
	\caption{Some products of pieces}
\end{figure}

The reader may wish to draw some of the remaining subdiagrams in order to 
familiarize himself or herself with their appearance.

Let us now assume that $x_{k-1} \cdot x_{k-2}$ also sees ILTC, and paste 
to $\sv'(k,k-1)$ the piece corresponding to $x_{k-2}$.  We observe that 
aside from the letters whose equality is forced in order to produce the 
ILTC, there can be no equality between the remaining letters in $x_k$, 
$x_{k-1}$, and $x_{k-2}$.  For instance, in multiplying $O \cdot Os'_{k-1} 
\cdot E$, $s'_{k-2}=s''_{k-1}$ cannot be, for otherwise 
$\{[s'_ks''_{k-1}],[s''_{k-1}s'_{k-1}],[s'_{k-1}s'_k]\}$ would be a 
circuit of length 3 containing two edges corresponding to edges of $C$, 
contradicting $C$'s achordality.

We make another observation.  Consider the portion of $\sv'(k,k-2)$ which 
lies between $s'_k$ and $s'_{k-2}=s''_{k-2}$; that is, the largest 
subdiagram of $\sv'(k,k-2)$ all of whose vertices lie on some simple path 
from $s'_k$ to $s'_{k-2}$.  This subdiagram has diameter 2: any two 
vertices in this subdiagram can be connected by a path $P$ of length at 
most 2.  Furthermore, the edges in $P$ can be chosen so that at most one 
edge in $P$ lies in $\sv'(k,k-1)$ and not $\sv'(k-1,k-2)$, and at most 
one edge lies in $\sv'(k-1,k-2)$ and not $\sv'(k,k-1)$.

We may generalize these observations as we continue to paste together the 
pieces corresponding to $x_k,x_{k-1},...,x_j$, as long as each word $x_i$ 
has no more than one long term and as long as $x_i \cdot x_{i-1}$ 
witnesses some non-exceptional case of ILTC.  We have the following lemma.

\begin{lemma} \label{noloops}
Suppose that each of the words $x_k,x_{k-1},...x_j$ has at most one 
long term, that each word $x_i$ sees non-exceptional ILTC with the 
following word, and that $j>1$.  Denote by $\sv''$ the portion of the 
subdiagram $\sv'(k,j)$ which lies between $s'_k$ and $s''_j$ in the 
sense described above.  Then the following all hold.

\vskip 2mm

\noindent 1. \ For every $i$, $j \leq i \leq k-1$, $s'_i$ and $s''_i$ lie 
in $\sv''$.

\vskip 2mm

\noindent 2. \ Besides the letters forced to be equal by ILTC (and, of 
course, $s'_i=s''_i$ in case $x_i$ is of type E or B), there is no 
equality between any of the letters in $\sv''$.

\vskip 2mm

\noindent 3. \ The diameter of $\sv''$ is at most $k-j$, and moreover any 
two vertices in $\sv''$ can be connected by a simple path $P$ so that for 
every subdiagram $\sv'(i,i-1)$ of $\sv''$ ($j+1 \leq i \leq k$), at most 
one edge of $P$ lies in $\sv'(i,i-1)$ and in no other such subdiagram.
\end{lemma}

\begin{proof}
We sketch a proof and leave the details to the reader.

The proof is essentially an induction on $k-j$.  In case $k-j=1$, 
the definition of ILTC and a glance at all possible pieces yields the 
desired conclusions.  We then assume that the result holds for all values 
of $k-j$ less than a given value $m$, and consider $x_kx_{k-1} \cdots 
x_{j+1} \cdot x_j$ where $m=k-j$.

(1) is easily proven by examining all possible pieces for products 
$x_{j+1} \cdot x_j$ and appealing to the inductive hypothesis regarding 
$\sv'(k,j+1)$.

(3) will also follow from the corresponding fact about $\sv'(k,j+1)$ once 
(2) is established for $\sv'(k,j)$.  For purposes of illustration, let 
$x_{j+1} \cdot x_j$ have the form $O \cdot O$ (other cases can be handled 
in a similar fashion).  Because ILTC occurs, $s''_{j+1}=s'_j$.  Assume to 
derive a contradiction that $s''_j$ is equal to some generator which lies 
in the portion of $\sv'(k,j+1)$ between $s'_k$ and $s''_{j+1}$.

First assume that $s''_j$ is not adjacent to $s'_j \in \sv'(k,j+1)$.  
Then by inductive hypothesis, the edge $[s'_js''_j]$ may be concatenated 
with a path $P$ of length at most $k-(j+1)=m-1$ to obtain a circuit $D'$, 
whose length is at most $m \leq k-1$.

Therefore we may inductively apply Theorem~\ref{cmatch} and find a circuit 
$D$ in $\sv$ to which $D'$ corresponds, edge for edge.  However, because 
$P$ can be chosen as in (3), either the edges of $\sv$ which correspond to 
those of $D'$ do not form a circuit at all, or they contradict the 
achordality of $C$.

The case in which $s''_j$ and $s'_j$ are adjacent in $\sv'(k,j+1)$ 
requires a different argument.  However, using the facts that the edges of 
$C$ must correspond as in Theorem~\ref{cmatch} to distinct edges of $C'$ 
and that $C$ is achordal, one may complete the proof in this case as well.  
The reader is encouraged to experiment with a few different cases in order 
to verify this claim.  (Compare the example for $j=k-2$ given before the 
statement of Lemma~\ref{noloops}.)
\end{proof}

Notice that we have assumed $j>1$ so that we can appeal to inductive 
hypothesis regarding the length of $C$.  If $j=1$ the argument in 
Lemma~\ref{noloops} will still go through as long as the equality that 
occurs is not $s'_k=s''_1$ or $s'_k=s'_1$.  Unfortunately, we must 
consider this case: suppose that $s''_i=s'_{i-1}$ for $i=1,...,k$.  Then 
if $x_i=s'_is''_i$ for every $i$, $x_kx_{k-1} \cdots x_1 = 1$ clearly 
holds, and there is ILTC between any two adjacent terms in this product!  
(This case is analogous to the case of the affine Euclidean groups 
that arose when $k=3$.)  However, we claim that in this case, we can 
either appeal to Lemma~\ref{starpiecing} or derive a contradiction.

The key here is that $x_kx_{k-1} \cdots x_j = s'_ks''_j$ for every 
$j=1,...,k$.  That is, in the ILTC that occurs in this product, every 
letter of each term $x_{k-1},x_{k-2},...,x_{j+1}$ is canceled, because 
every product $x_i \cdot x_{i-1}$ is of the form $O \cdot O$ and each term 
$x_i$ has length 2.  The number of letters in this reduced product which 
follow $s'_k$ never increases.

Let every product $x_i \cdot x_{i-1}$ be of a type which permits at most 
one pair of letters to cancel (i.e., not of types 4 or 5).  Assuming that 
the number of letters following the last occurrence of $s'_k$ never 
increases (as above), none of the words $x_i$ can contain an even long 
term.  Moreover, we see also that there can be no products of the type $O 
s'_i \cdot O$, and $s'_is''_i$ cannot have order greater than 3.

Even if we allow products $x_i \cdot x_{i-1}$ which admit more than one 
pair of letters to cancel, the length of the product $x_k \cdots x_j$ will 
increase.  Indeed, even in cases 4 and 5, cancellation of more than one 
pair of letters implies (by the two-dimensionality of $W$) that one of the 
long terms involved in the cancellation is indeed quite long.  For 
example, suppose that $x_i \cdot x_{i-1}$ is of type $O \alpha 
\cdot E$, so $x_i \cdot x_{i-1} = v_{s'_is''_i} \alpha \cdot u_{\beta 
s'_{i-1}}$, where $s'_is''_i$ has order $m$ and $\beta s'_{i-1}$ has order 
$n$.  If more than a pair of letters cancels, $\alpha=\beta$ and 
$s''_i=s''_{i-1}$.  Then $n=4 \Rightarrow m \geq 7$ and $m=3 \Rightarrow n 
\geq 6$.  Therefore none of these terms can occur if $x_k \cdots x_1 = 1$ 
holds.

We are therefore left with only a few types of product that $x_i \cdot 
x_{i-1}$ can be, all involving only terms of types $O$ and $B$.  First 
assume that there are no blank terms.  Arguing as in the proof of 
Lemma~\ref{noloops}, we can show that each subdiagram $\sv'(k,j)$ appears 
as a sequence of triangles, each sharing an edge with the last.  As in 
that proof, once we have three such triangles in a row, we can apply 
Lemma~\ref{starpiecing} and conclude.

Note that there cannot be two non-trivial blank terms in a row, for if 
there were, they would not cancel with each other, since we are assuming 
that there is no short term cancellation.  Each blank term $x_i$ which 
occurs between two odd terms $x_{i+1}$ and $x_{i-1}$ must admit 
cancellation with both $x_{i+1}$ and $x_{i-1}$.  Moreover, as above, one 
can rule out products $O s'_i \cdot B$ and $O \alpha \cdot B$ (where 
$\alpha s'_i = s'_i \alpha$).  Thus the above blank term $x_i$ must be 
$s''_{i+1}s'_{i-1}$.  If $x_i$ is a non-trivial blank term following a 
string of trivial blank terms, it must be followed by an odd term 
$x_{i-1}$, so that $x_i=s'_{i-1}$; similarly, if $x_i$ is a 
non-trivial blank term which precedes a string of trivial blank 
terms, $x_i=s''_{i+1}$.  Finally, there can be no more than two terms of 
type $O$ in a row, as otherwise, we would obtain three consecutive 
triangles to which Lemma~\ref{starpiecing} could be applied, as in the 
previous paragraph.

Using this information one can piece together the subdiagram 
corresponding to the product $x_k \cdots x_1$; it has a rather regular 
form, consisting of a sequence of edges, triangles, and ``diamonds'' 
(pairs of triangles sharing an edge), each such component sharing a single 
vertex with the last.  In any such configuration, the assumption that 
$s''_1=s'_k$ will yield a shorter achordal circuit $D'$ in $\sv'$ to which 
Theorem~\ref{cmatch} can be applied.  This circuit must correspond with a 
circuit $D$ in $\sv$.  However, as in the proof of Lemma~\ref{noloops}, 
the edges of $\sv$ which correspond to those of $D'$ either do not form a 
circuit at all or contradict the achordality of $D$.  If we assume that 
$s'_1=s'_k$ instead, one can see that $s'_1=s''_1$ must hold as well, for 
given the above form of each $x_i$, $x_k \cdots x_1=1$ would not hold if 
$s'_1 \neq s''_1$.  In this case, we may again construct a shorter 
achordal circuit $D'$ to which Theorem~\ref{cmatch} can be applied 
inductively to obtain a contradiction.

Thus if $j=1$ in the above chain of ILTC, we have proven 
Proposition~\ref{noiltc}.  If $j \neq 1$, we have computed a reduced word 
which forms a ``prefix'' for $x_k \cdots x_1$.  We will now argue that in 
completing the product $x_k \cdots x_j \cdot x_{j-1} \cdots x_1$, there 
can be almost no cancellation of this prefix.  Essentially, we show how 
we can continue to multiply words $x_i$ (possibly with a great deal of 
cancellation) until a new ILTC chain as above is encountered, and then 
repeat the process until $x_k \cdots x_1$ is obtained.

If the ILTC chain above ends at $x_j$ ($j \geq 2$), it does so for one of 
the following reasons.

\vskip 2mm

1. \ $x_j \cdot x_{j-1}$ is already reduced (there is no cancellation),

\vskip 2mm

2. \ $x_j \cdot x_{j-1}$ sees exceptional ILTC and $x_{j-1}$ is completely 
canceled in this product,

\vskip 2mm

3. \ $x_j \cdot x_{j-1}$ sees ILTC and $x_{j-1}$ has at least two long 
terms.

\vskip 2mm

Consider for a moment the second case.  In this case, as mentioned 
immediately following the definition of ILTC, it is easily shown that 
there is no further cancellation between $x_k \cdots x_{j-1}$ and 
$x_{j-2}$ once $x_{j-1}$ is completely canceled.  Therefore, this case can 
be argued in much the same way as the first, with $x_{j-2}$ in place of 
$x_{j-1}$.

Suppose first that $x_k \cdots x_j \cdot x_{j-1}$ admits no further 
cancellation.  If $j=2$ we are done.  Otherwise, we consider the product 
$x_{j-1} \cdot x_{j-2}$.  If there is no cancellation here either, we may 
continue by considering $x_{j-2} \cdot x_{j-3}$, and now we have the 
additional advantage of $x_{j-1}$ serving as a ``buffer'' between $x_j$ 
and $x_{j-2}$ which effectively forbids further cancellation of letters in 
$x_k \cdots x_j$.  If $x_{j-1} \cdot x_{j-2}$ admits ILTC instead, either 
we begin a new chain of ILTC between words $x_i$ with at most one long 
term each or $x_{j-2}$ has at least two long terms.  In the latter case, 
the first long term in $x_{j-2}$ again serves as a buffer preventing 
further cancellation with $x_k \cdots x_j$ when further terms $x_{j-3}$, 
$x_{j-4},...$ are multiplied.

Finally, it is possible that $x_{j-1} \cdot x_{j-2}$ admits CLTC.  If 
$x_{j-1}$ has more than 2 terms, the first terms serve as a buffer, as 
above, preventing further cancellation with $x_k \cdots x_j$.  Otherwise, 
we must be more careful.

The following lemmas are useful when considering CLTC between terse words.

\begin{lemma} \label{CLTCnotmuch}
Suppose that $x_i \cdot x_{i-1}$ admits CLTC.  Then in multiplying $(x_k 
\cdots x_i) \cdot (x_{i-1} \cdots x_1)$, at most one one pair of long 
terms (those admiting the CLTC) and one additional pair of letters 
cancels.  (The letters from this additional pair may be in either 
$x_{i+1}$ or $x_{i-2}$.)
\end{lemma}

Lemma~\ref{CLTCnotmuch} essentially tells us that not much more than the 
completely canceled long terms cancels.  It is proven by straightforward 
arguments similar to those used to prove Proposition~\ref{xicancel}.

\begin{lemma} \label{noloops2}
Suppose that each of the words $x_l,x_{l-1},...,x_r$ is terse, that each 
word $x_i$ ($l \leq i \leq r$) sees CLTC with the following word, and that 
$r>1$.  Then the vertices $s'_l,s'_{l-1},...,s'_r$ and 
$s''_l,s''_{l-1},...s''_r$ all lie on a simple path $P$ of length at most 
$r-l+2$, and equality between any of these letters only occurs when forced 
by CLTC.  Moreover, the edges 
$[s'_ls''_{l-1}],[s'_{l-1}s''_{l-2}],...,[s'_{r+1}s''_r]$ lie in this 
order on $P$, with no intervening edges.
\end{lemma}

Lemma~\ref{noloops2} is proven in much the same way as was 
Lemma~\ref{noloops}.  Analogously, it allows us to multiply successive 
terse words $x_i,x_{i-1},...$ as long as each such word admits CLTC with 
the next.  As was the case with Lemma~\ref{noloops}, the subdiagram of 
$\sv'$ that emerges from Lemma~\ref{noloops2} closely resembles the 
circuit $C'$ we wish to construct.

The essence of the remaining proof is as follows.  Having finished with 
the initial chain $x_k \cdots x_j$ of words exhibiting ILTC, either we 
begin a new chain of ILTC or we begin a chain of CLTC, perhaps with a 
buffer in between (provided by non-terse words or words $x_i$ and 
$x_{i-1}$ which do not cancel).  After completing the next chain of ILTC 
or CLTC, the same occurs, and we repeat the process.

Arguments similar to those used to prove Lemmas~\ref{noloops} and 
\ref{noloops2} show that in transitioning from one chain of cancellation 
to another {\it without} a buffer in between, we maintain the 
``diametric'' property described by those lemmas.  That is, the diameter 
of the portion of the subdiagram $\sv'(k,j)$ lying between $s'_k$ and 
$s''_j$ is small enough to allow an inductive appeal to 
Theorem~\ref{cmatch} which forbids any backtracking produced by equality 
of vertices of $\sv'(k,j)$ other than that forced by ILTC or CLTC.  On the 
other hand, if a buffer {\it does} appear between two chains, this buffer 
prevents us from canceling every letter of the preceding chain, so that 
ultimately the product $x_k \cdots x_1$ cannot be trivial.

This discussion has focussed upon the first case mentioned above (in which 
$x_j~\cdot~x_{j-1}$ admits no cancellation).  Clearly the third case can 
be handled in a similar fashion.

There is one difficulty which must be overcome, and it concerns the single 
pair of letters that could be canceled in addition to the CLTC in 
Lemma~\ref{CLTCnotmuch}.  If this additional cancellation occurs at the 
end of the initial chain of ILTC, such cancellation may ``expose'' the 
letter $s'_k$ at the end of $x_k \cdots x_{j-1}$, and this letter could 
then be canceled when $x_1$ is at last multiplied with $x_k \cdots x_2$.

In most such cases, we may solve the problem by an appeal to the inductive 
hypothesis of either Theorem~\ref{cmatch} or Lemma~\ref{starpiecing}.  For 
example, consider the following case, in which $k=6$: $x_6 = s'_6s''_6 = 
s'_6s'_5$, $x_5 = s'_5s''_5 (\alpha s'_5)^3$, $x_4 = (\beta s'_4)^3 
s''_5$, $x_3 = x_2 = 1$, $x_1 = s'_6$.

Here $x_6 \cdots x_1 = 1$, given that $\alpha=s'_4$, $\beta=s'_5$, and 
$s'_6s'_1=s'_1s'_6$.  However, one can draw the subdiagram of $\sv'$ 
determined by the above equalities and see that there is a circuit of 
length 5 to which we may apply Theorem~\ref{cmatch} inductively in order 
to contradict the achordality of $C$.  (Alternatively, one could apply 
Lemma~\ref{starpiecing} in this case as well.)

Schematically, the example above begins with the short ILTC chain $x_6 
\cdot x_5$ of type $O \cdot OEs'_5$.  There are a number of other cases 
where the initial ILTC chain $x_k \cdots x_j$ yields a prefix ending with 
$s'_k$, and in almost all of these cases one can apply at least one of the 
two arguments mentioned above.  As when $x_k \cdots x_1$ was a single ILTC 
chain, a problem arises when there are not enough adjacent triangles to 
apply Lemma~\ref{starpiecing} inductively, and when there is no place to 
apply Theorem~\ref{cmatch} to a shorter circuit.  This situation is 
illustrated by the following example.

Let $x_6 = s'_6s''_6$, $x_5 = s''_6$, $x_4 = (s'_4s''_4)^3s'_4$, $x_3 = 
(s'_3s''_3)^3s'_3$, $x_2=1$, and $x_1 = s'_6$, where $s'_1 \neq s'_6$ and 
$s'_5 \neq s''_6$ are commuting pairs, $s'_3=s'_4$, and $s''_3=s''_4$.  
Then the subdiagram of $\sv'$ corresponding to the product $x_6 \cdots 
x_1$ contains only a single pair of triangles which meet at an edge 
(ruling out use of Lemma~\ref{starpiecing}) and no circuit of length less 
than 6 to which an application of Theorem~\ref{cmatch2} is helpful.

We assert that in any case to which we can apply neither 
Lemma~\ref{starpiecing} nor Theorem~\ref{cmatch}, a similar configuration 
arises in $\sv$ and in $\sv'$.  Namely, there exists a vertex $\gamma$ in 
$\sv$ such that for some $i$, the vertices $s_i$, $s_{i-1}$, $s_{i-2}$, 
and $\gamma$ appear as in Figure~3.a, where both $[s_is_{i-1}]$ and 
$[s_{i-1}s_{i-2}]$ have labels greater than 2, and $[s_{i-1}\gamma]$ has 
an odd label.  Finally, the two triangles shown in Figure~3.a correspond 
as in Theorem~\ref{cmatch} to the triangles in $\sv'$ shown in Figure~3.b.
We can argue much as we did when $x_k \cdots x_1$ was a single ILTC chain 
to prove that this configuration must arise.

\begin{figure}[hbt]
	\begin{center}
		\input{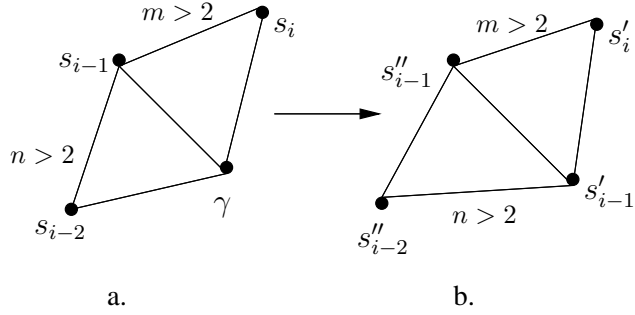}
	\end{center}
	\caption{A final case}
\end{figure}

Now let $\bar{w}_i$ conjugate each vertex in $\{s_i,s_{i-1},\gamma\}$ to 
the appropriate vertex in $\{s'_i,s'_{i-1},s''_{i-1}\}$ and 
$\bar{w}_{i-1}$ conjugate each vertex in $\{s_{i-1},s_{i-2},\gamma\}$ to 
the appropriate vertex in $\{s'_{i-1},s''_{i-1},s''_{i-2}\}$.  (These are 
the words whose existence is guaranteed by Theorem~\ref{cmatch2} in the 
case $k=3$.)  Let $w_i$ and $w_{i-1}$, as always, denote the elements 
conjugating the edges $[s_is_{i-1}]$ and $[s_{i-1}s_{i-2}]$, so that 
$x_{i-1}=w_iw_{i-1}^{-1}$.

\begin{lemma} \label{finalcase}
Let the configuration pictured in Figure~3 appear.  Then, using the 
notation from the preceding paragraph, $x_{i-1} \neq 
v_{s'_{i-1}s''_{i-1}}$.
\end{lemma}

In the above example, $i=1$, and Lemma~\ref{finalcase} finishes the proof 
of Proposition~\ref{noiltc} by contradicting the fact that 
$x_k=v_{s'_6s''_6}$.

\begin{proof}

Assume that $s'_is''_{i-1}$ and $s'_{i-1}s''_{i-2}$ both have odd order.
Using Lemma~\ref{edgenormalizers}, it is easy to prove that $w_i 
\bar{w}_i^{-1} \in \{1,v_{s'_is''_{i-1}}s'_i\}$, 
$w_{i-1}\bar{w}_{i-1}^{-1} \in \{1,v_{s'_{i-1}s''_{i-2}}s'_{i-1}\}$, and 
$\bar{w}_i\bar{w}_{i-1}^{-1}=v_{s'_{i-1}s''_{i-1}}s'_{i-1}$.

Therefore

$$x_{i-1} = w_iw_{i-1}^{-1} = w_i \bar{w}_i^{-1} \bar{w}_i 
\bar{w}_{i-1}^{-1} \bar{w}_{i-1} w_{i-1}^{-1},$$

\noindent which from the above computations can be one of four different 
products, none of which is $v_{s'_{i-1}s''_{i-1}}$.

\end{proof}

This completes the proof of Proposition~\ref{noiltc} by taking care of the 
last remaining cases.

We are now close to a proof of Theorem~\ref{cmatch}.

Recall we have assumed throughout $z_{i-1,i}=1$ for all $i$.  We note now 
that the arguments we have developed above prove that $z_{i-1,i}=1$ must 
hold.  Suppose that $x_i = x'_i z_{i-i,i}^{-1}$, where $z_{i-1,i} \neq 1$.  
Then the product $x_i \cdot x_{i-1}$ can be treated much like a case of 
ILTC (as indeed $x_i \cdot x_{i-1}$ can never be trivial if $z_{i-1,i} 
\neq 1$), and we can modify all of the arguments above to take this 
possibility into account.

We may also use the arguments above to show that every word $x_i$ is 
terse; otherwise some $x_i$ would serve as a ``buffer'' which would 
prohibit the product $x_k \cdots x_1$ from being trivial.

The remaining words must all exhibit CLTC with one another.  In order to 
avoid contradicting Lemma~\ref{noloops2}, there can be no equality (except 
that forced by CLTC) between any elements $s'_k,s''_k,...,s'_2,s''_2$; at 
the last step, in multiplying $x_k \cdots x_2 \cdot x_1$, we must complete 
the circuit $C'$ of length $k$, corresponding to $C$.  Therefore 
Theorem~\ref{cmatch} is proven in case there is no cancellation of short 
terms.

However, we can modify all of our arguments to take care of such 
cancellation as well.  If $x_i=x'_i \alpha_1 \cdots \alpha_l$ and 
$x_{i-1} = \alpha_l \cdots \alpha_1 x'_{i-1}$ where each $\alpha_j$ is a 
short term, then we can cancel all of the letters $\alpha_j$ and apply our 
ILTC and CLTC arguments to $x'_i$ and $x'_{i-1}$ instead.  Intermediate 
short term cancellation that arises after long terms have been canceled is 
met in a similar fashion.

To prove Theorem~\ref{cmatch2}, we observe that in performing complete 
cancellation of words $x_i$ and $x_{i-1}$, the precise form of the long 
terms in each of these words is forced, as in Section 4 when all cases in 
which $k \in \{3,4\}$ were considered.  Also, Lemmas~\ref{whichvertex} 
and \ref{overthehump} were proven in complete generality, and as indicated 
in that section, Lemmas~\ref{whichword}, \ref{lotsocancellation}, and 
\ref{lotsocancellation2} follow from arguments similar to the proof of 
Proposition~\ref{noiltc} in case $k \geq 5$.  Applying these results 
concludes our proof.

\section{Piecing circuits together}

We have now established Theorem~\ref{cmatch} and Theorem~\ref{cmatch2}.  

Assume now that $\sv$ is odd-edge-connected, and let $C_1$ and $C_2$ be 
two circuits in $\sv$ which share at least one edge.  Let $C'_1$ and 
$C'_2$ be the circuits in $\sv'$ to which $C_1$ and $C_2$ correspond as in 
Theorem~\ref{cmatch}, respectively.  Then $C'_1$ and $C'_2$ share edges 
corresponding to the common edges of $C_1$ and $C_2$.  If there is more 
than one common edge, there is no ``twisting'' at any edge, so $C_1 \cup 
C_2$ and $C'_1 \cup C'_2$ are isomorphic as edge-labeled graphs.

Suppose there is a single common edge, $[st]$, with an odd label (if 
the label of this edge is even, there can be no twisting at this edge).  
Because $\sv$ is odd-edge-connected, the removal of this edge does not 
disconnect the diagram $\sv$, and we can find a simple path $P$ in $\sv$ 
leading from a vertex in $x_1 \in C_1 \setminus \{s,t\}$ to a vertex $x_2 
\in C_2 \setminus \{s,t\}$.  Moreover, we can choose this path so that 
$x_1$ and $x_2$ are as close to $s$ as possible, and so that $P$ is of 
minimal length among paths satisfying this first condition.  (Both of 
these conditions can be met by replacing subpaths of $P$ with shorter 
paths as needed.)  Denote by $P_1$ the path from $x_1$ to $s$, and by 
$P_2$ the path from $s$ to $x_2$.  Then the path $P_1P_2P^{-1}$ is a 
circuit, $D$, and $D$ is achordal, except perhaps for edges $[sy_i]$, 
where $y_i \in P$.  Subdivide $D$ by adding these edges, as necessary, 
into circuits $D_1,...D_r$.

We have obtained a picture very similar to that considered in 
Lemma~\ref{starpiecing}.  An argument almost identical to the proof of 
that lemma now shows that twisting can occur neither at $[sy_i]$ for any 
$i$, nor at $[st]$.  Therefore, $C_1 \cup C_2$ and $C'_1 \cup C'_2$ are 
isomorphic as edge-labeled graphs.

If $\sv$ is odd-edge-connected, it is easy to see that every vertex lies 
on some achordal circuit.  Thus by piecing together achordal circuits 
which share at least one edge in the manner described above, we prove 
Theorem~\ref{edgeconnected}.  If $\sv$ is still one-connected but the 
removal of some odd edge $e$ disconnects $\sv$, we can induct on the 
number of ``odd-edge-indecomposable'' pieces into which $\sv$ may be 
divided by removing such edges in order to prove Theorem~\ref{main}.

Now suppose that $\sv$ is connected but not one-connected.  In this case, 
we can apply the same technique as used by M\"uhlherr and Weidmann in 
\cite{MuWe} to prove their Main Theorem.  (See Section 8 of \cite{MuWe}.  
Theorem~\ref{edgeconnected} here serves as the base case for the inductive 
proof.)  This technique draws heavily upon the results of \cite{Deo} and 
\cite{MiTs}.  The latter paper details a canonical decomposition for a 
given Coxeter group, arrived at through an application of Bass-Serre 
theory.  As was done in \cite{MuWe}, we may prove Theorem~\ref{main} by 
inducting upon the number of ``vertex-indecomposable'' pieces into which 
the diagram $\sv$ can be broken by removing separating vertices.

Appealing to \cite{FR1} and \cite{FR2} (as was done in \cite{MuWe}), we 
complete the proof in case $\sv$ is not connected.  This concludes the 
verification of Theorem~\ref{main}.

There is an immediate corollary of Theorem~\ref{main}, regarding the 
structure of Artin groups.  Recall that the Artin group $A(S)$ 
corresponding to a given Coxeter presentation $W \cong \langle S \ 
| \ R \rangle$ is found merely by deleting from $R$ the relators $s^2$, $s 
\in S$.  (Therefore there is an epimorphism from $A(S)$ to $W$ which maps 
each element $s^2$ to the identity, for $s \in S$.)  Clearly the diagram 
$\sv$ for $(W,S)$ completely determines the group $A(S)$ as well as the 
group $W$.  We may define reflections, rigidity, and reflection rigidity, 
in Artin groups in exactly the same way we have defined them for Coxeter 
groups.

From Theorem~7.2 of \cite{BMMN}, we derive the following result.

\begin{theorem} \label{artin}
Let $A(S)$ be the Artin group corresponding to the two-dimensional Coxeter 
system $(W,S)$, with diagram $\sv$.  Then $A(S)$ is reflection rigid, up 
to diagram twisting.  (That is, given any other Coxeter system $(W,S')$ 
such that $A(S)$ and $A(S')$ yield the same reflections, the diagram for 
$A(S')$ can be derived from $\sv$ by a sequence of diagram twists.)
\end{theorem}

What else can be said?  As we have seen, achordal circuits in $\sv$ are 
nearly strongly rigid; conjugating words for the various vertices in such 
a circuit $C$ differ only by products $\alpha_1 \cdots \alpha_l$, for 
appropriately chosen $\alpha_i$.  If $\sv$ has no edges labeled 2, every 
achordal circuit is strongly rigid, and arguing as in \cite{Ba3}, we 
recover another fact proven in \cite{MuWe}: if $W$ is a skew-angled 
reflection independent Coxeter group and the diagram $\sv$ for the system 
$(W,S)$ is edge-connected, then $W$ is strongly rigid.

We can still say something when $\sv$ contains edges labeled 2.  Let 
$(W,S)$ be an arbitrary Coxeter system, with diagram $\sv$, and let $s$ 
be a vertex in $\sv$.  As in \cite{Ba2}, \cite{Ba3}, and \cite{BaMi1} we 
define the {\it 2-star}, ${\rm st}_2(s)$, of $s$ to be the set of vertices

$$\{s\} \cup \{t \in \sv \ | \ [st] \ {\rm is \ an \ edge \ labeled} \ 2 
\} \subseteq \sv.$$

We have the following theorem.  (Compare this with the main theorem of 
\cite{Ba3}).

\begin{theorem} \label{sr}
Suppose that $(W,S)$ is a reflection independent two-dimensional Coxeter 
system whose diagram $\sv$ has at least 3 vertices.  Then $\sv$ is 
strongly rigid if $\sv$ is edge-connected and there are no vertices 
$s,t_1,t_2 \in \sv$ such that the removal of ${\rm st}_2(s)$ separates 
$\sv$ into at least 2 components, $t_1,t_2 \not\in {\rm st}_2(s)$, and 
$t_1$ and $t_2$ lie in different components of the full subdiagram of 
$\sv$ induced by the vertices $\sv \setminus {\rm st}_2(s)$.
\end{theorem}

\begin{proof}
We have already seen that if $\sv$ has three 3 vertices, $W$ is strongly 
rigid.  Therefore we may assume that $\sv$ has at least 4 vertices.

Let $\sv$ satisfy both of the conditions put forth in the statement of 
the theorem, and let $(W,S')$ be another Coxeter system for $W$, with 
diagram $\sv'$.  Note that because the diagram is edge-connected (and 
therefore odd-edge-connected), Theorem~\ref{edgeconnected} shows that it 
is reflection rigid, and therefore rigid, because $W$ is assumed to be 
reflection independent.  Therefore $\sv$ and $\sv'$ are isomorphic, and 
the achordal circuits in these diagrams match up as in 
Theorem~\ref{cmatch2}.

We first claim that every edge (and therefore every vertex) of $\sv$ must 
lie on an achordal circuit.  In fact, it is easy to show that if $[st]$ 
did not lie on any circuit, then removing $[st]$ from $\sv$ would separate 
the diagram, contradicting our hypotheses.  Thus every edge lies on a 
circuit, which can be shortened, if needed, to make it simple and 
achordal.  We claim that to a given achordal circuit $C$ there is an 
element $w_C \in W$ which conjugates each vertex of $C$ to the appropriate 
vertex of $\sv'$.

From Section 4, this is so if the circuit has length 3.  Thus, we may 
assume that $C = \{[s_1s_2],...,[s_ks_1]\}$, $k \geq 4$.  Let $s_i \neq 
s_j$ be vertices on $C$.  Let $w_i$ and $w_j$ be elements which conjugate 
$s_i$ to $\hat{s}_i$ and $s_j$ to $\hat{s}_j$, respectively.  We assume 
that $i>j$.  As in Theorem~\ref{cmatch2}, $w_iw_j^{-1}$ is a product of 
elements $\alpha_i \in \sv'$ which either satisfy that theorem's 
separation condition or lie in $\{\hat{s}_i,\hat{s}_j\}$.

In case $\sv = C$, the result follows from \cite{ChDa}, so we may suppose 
that $\sv \neq C$.  Consider a letter $\alpha$ appearing in the ratio 
$w_i^{-1}w_j$ which does not lie on $C'$.  The removal of ${\rm 
st}_2(\alpha)$ separates $C'$ into various subarcs, but because the second 
condition on $\sv$ in the statement of the theorem does not obtain, 
removing ${\rm st}_2(\alpha)$ does not disconnect $\sv'$.  By 
two-dimensionality, $\alpha$ is not adjacent to at least two letters in 
$C'$ by edges labeled 2, so $C' \setminus {\rm st}_2(\alpha)$ is not 
empty.  Therefore given any two subarcs $P_1$ and $P_2$ into which $C'$ is 
divided by removing ${\rm st}_2(\alpha)$, there is a path $P$ lying 
completely in $\sv' \setminus {\rm st}_2(\alpha)$ which connects $P_1$ and 
$P_2$.  Let $P_1$ and $P_2$ be ``adjacent'' subarcs of $C'$, lying on 
either side of the vertex $\hat{s}_i$ (where $[\hat{s}_i\alpha]$ is an 
edge labeled 2).  By replacing portions of $P$ with the appropriate paths, 
we may assume that the endpoints of $P$ are as close as possible to 
$\hat{s}_i$ and that $P$ is as short as possible among all paths lying in 
$\sv' \setminus {\rm st}_2(\alpha)$ connecting these endpoints.  Fix such 
a path $P$, with endpoints $x_1 \in P_1$ and $x_2 \in P_2$.  Denote by 
$Q_l$ the subpath of $P_l$ from $x_l$ to $\hat{s}_i$, for $l=1,2$.  Then 
by our choice of $P$, the circuit $D=PQ_2Q_1^{-1}$ is achordal, except 
perhaps for edges $[\hat{s}_iy]$, for some $y \in P$.  Add such edges as 
needed to subdivide $D$ into achordal circuits $D_1,...D_r$.  An argument 
we have now seen twice before implies that since $\hat{s}_i$ is the only 
vertex in any of these circuits which is adjacent to $\alpha$ by an edge 
labeled 2, $\alpha$ cannot appear in any of the ratios of conjugating 
elements for the edges in each circuit $D_l$.  Therefore if 
$\hat{s}_{i_1}$ lies in $P_1$ and $\hat{s}_{i_2}$ lies in $P_2$, the ratio 
of the conjugating elements associated to these vertices cannot contain 
$\alpha$.

Repeating this procedure for every pair of adjacent subarcs $P_1$ and 
$P_2$, and then for every element $\alpha$ whose 2-star separates $C'$, we 
see that no such $\alpha$ can occur.  Therefore the only letters $\alpha$ 
that appear in $w_iw_j^{-1}$ are $\hat{s}_i$ and $\hat{s}_j$.  Thus we 
have effectively reduced the problem to the case in which $\sv=C$.

Now consider two achordal circuits $C_1$ and $C_2$ which share at least 
one edge.  Using arguments almost entirely like those just applied, one 
can show that $C_1$ and $C_2$ share a common conjugating element.  
Therefore, since every vertex in $\sv$ lies on some achordal circuit, we 
have proven the theorem.
\end{proof}

\end{document}